\documentclass[12pt]{article}
\usepackage{amsfonts}
\usepackage{amsmath,amssymb} 
\newtheorem{thm}{Theorem}[section]

\newtheorem{cor}[thm]{Corollary}
\newtheorem{lem}[thm]{Lemma}
\newtheorem{propn}[thm]{Proposition}
\newtheorem{rem}[thm]{Remark}

\def\R{{\mathbb R}}
\def\E{{{\mathbb E}\,}}
\def\N{{\mathbb N}}
\def\P{{\mathbb P}}
\def\Z{{\mathbb Z}}
\def\sL{{\cal L}}
\def\calam{{\cal A}^{\lambda}}
\def\sA{{\cal A}}
\def\sS{{\cal S}}
\def\sP{{\cal P}}
\def\ce{{\cal E}}
\def\cf{{\cal F}}
\def\cs{{\cal S}}

\def\proof{{\noindent{\sc Proof.}~}}

\def\qed{{\hfill $\square$ \bigskip}}

\def\eps{\varepsilon}
\def\vp{\varphi}
\def\bp{{\bar p}}

\def\fal{{~~~\mbox{for all }~}}
\def\tfrac#1#2{{\textstyle {#1\over #2}}}
\def\ol{\overline}
\def\sgn{{\mbox{sgn}\,}}

\makeatletter
\@addtoreset{equation}{section}

\makeatother
\begin{document}
\title{Symmetric Markov chains on $\Z^d$
with unbounded range}
\author{Richard F. Bass
\thanks{Research partially supported by NSF grant DMS0244737.}
\qquad and \qquad
Takashi Kumagai
\thanks{Research partially supported by Ministry of Education, Japan,
Grant-in-Aid for Scientific Research for Young Scientists (B) 16740052.}}

\maketitle

\begin{abstract}
We consider symmetric Markov chains on $\Z^d$ where we do {\bf not}
assume that the conductance between two points must be zero if the points are
far apart. Under a uniform second moment condition on the conductances,
we obtain upper bounds on the transition probabilities, estimates for
exit time probabilities, and certain lower bounds on the transition probabilities.
We show that a uniform Harnack inequality holds if an additional assumption
is made, but that without this assumption such an inequality need not hold.
We establish a central limit theorem giving conditions for a sequence of 
normalized symmetric Markov chains
to converge to a diffusion on $\R^d$ corresponding to an elliptic
operator in divergence form.

\end{abstract}

\def\Var{\mathop {\rm Var\, }}
\def\half{{1/2}}

\section{Introduction}

Let $X_n$ be a symmetric Markov chain on $\Z^d$. We say that $X_n$ has {\it bounded}
range if there exists $K>0$ such that $\P(X_{n+1}=y \mid X_n=x)=0$ whenever
$| y-x|\geq K$. The range is {\it unbounded} if for every $K$ there exists $x$ and $y$
(depending on $K$) with $|x-y|>K$  such that $\P(X_{n+1}=y \mid X_n=x)>0$.
 There is a great deal known about Markov chains on graphs when
the chains have bounded range. The purpose of this paper is to
obtain results for Markov chains on $\Z^d$ that have  unbounded range.

Suppose $C_{xy}$ is the conductance between $x$ and $y$. We impose a condition on
$C_{xy}$ (see (A3) below) which essentially says that the $C_{xy}$ satisfy a uniform second
moment condition. 
Let $Y_t$ be the continuous time Markov chain on $\Z^d$ determined
by the $C_{xy}$, while $X_n$ is the discrete time Markov chain determined by these
conductances. The transition probabilities for the Markov
chain $X$ are defined by
\[
\P^x(X_1=y)=\frac{C_{xy}}{\sum_z C_{xz}},
\]
while the process $Y_t$ is the Markov chain that has the same jumps as $X$ but where the
times between jumps are
independent exponential random variables with parameter 1.
When (A3) holds, together with two very mild regularity conditions, 
we obtain upper bounds on the transition probabilities 
of the form
\[ \P(Y_t=y\mid Y_0 =x) \leq  ct^{-d/2} \]
 and some corresponding lower
bounds when $x$ and $y$ are not too far apart. Unlike the case of bounded range,
reasonable  universal bounds of Gaussian type need not hold when the range is unbounded.
We also obtain bounds on the exit 
probabilities $\P(\sup_{s\leq t} |Y_s-x|>\lambda t^{1/2})$.

We say a uniform Harnack inequality holds for $X$ if 
whenever $h$ is nonnegative and harmonic for the Markov chain $X$ in the ball 
$B(x_0,R)$ of radius $R>1$
about a point $x_0$, then
\[ h(x)\leq Ch(y), \qquad |x-x_0|, |y-x_0|<R/2, \]
where $C$ is independent of $R$.
Even when $X_n$ is a random walk, i.e., the increments $X_n-X_{n-1}$ form an
independent identically distributed sequence, a uniform Harnack inequality need
not hold. However, if we impose an additional strong assumption (see (A4)) on the conductances,
then we can prove such a Harnack inequality.

We prove that if we have   Markov chains $X^{(n)}$  on $\Z^d$
satisfying Assumption (A3) uniformly in $n$, 
the sequence of processes $X^{(n)}_t=X_{[nt]}/\sqrt n$ is 
tight in the space $D[0,\infty)$
of right continuous, left limit functions, and all subsequential limit points are 
continuous processes. 
Under an additional
condition on the conductances (A5) (different than the one needed for the Harnack
inequality), we then 
show that the $X^{(n)}_\cdot$ converge 
weakly as processes
to the law of the diffusion corresponding to an  elliptic operator 
\[
{{\cal L}} f(x)=\sum_{i,j=1}^d 
\frac{\partial}{\partial x_i} \Big( a_{ij}(\cdot) \frac{\partial f}{\partial x_j}(\cdot)\Big)
(x)
\]
in divergence form.
The exact statement is given by Theorem \ref{clt}.

In the case of bounded range Markov chains on $\Z^d$ some of our estimates 
have been obtained
by \cite{SZ}, and we obviously owe a debt to that paper. 
Not all of their methods extend to the unbounded case, however.
In particular, 
\begin{enumerate}
\item   New techniques were needed to obtain the exit probability estimates.
\item   A new method was needed to obtain lower bounds for the 
process killed on exiting a ball.
This method should
apply in many other instances, and is of interest in itself.
\item   Harnack inequalities in the case of unbounded range are quite a bit
more subtle, and this section is all new.
\item   In the proof of the central limit theorem, new methods
were needed to handle the case of unbounded range. 
Moreover, even in the bounded range
case our result covers more general situations.
\end{enumerate}

There are many versions of central limit theorems that investigate the asymptotic behavior of
$ \sum_{i=1}^n f(X_i) $
when $X_n$ is a symmetric Markov chain on a graph. These are quite
different from the central limit theorem of this paper.  
Our formulation has much more in common with the work of Stroock and Varadhan \cite{SV}, Chapter 11.
There they consider certain non-symmetric chains and show convergence to the law of
a diffusion corresponding to an operator in nondivergence form:
\[ \sL f(x)=\sum_{i,j=1}^d a_{ij}(x) \frac{\partial^2 f}{\partial x_i \partial x_j} (x)
+\sum_{i=1}^d b_i(x) \frac{\partial f}{\partial x_i}(x). \]
Our result is the analogue for symmetric chains and operators in divergence form.

The next section sets up the notation and framework and states the assumptions we need.
Section 3 has the exit time and hitting time estimates, Section 4 has the lower bounds, and
Section 5 discusses the Harnack inequality.  Our central limit
theorem is proved in Section 6.

The letter $c$ with or without subscripts and primes will denote
finite positive constants whose exact value is unimportant and
which may change from line to line.

\section{Framework}

We let $|\cdot|$
be the Euclidean norm and $B(x,r):=\{y\in \Z^d: |x-y|< r\}$. 
We sometimes write $|A|$ for the cardinality of a set $A\subset \Z^d$.

For each $x,y\in \Z^d$ with $x\ne y$, let $C_{xy}\in [0,\infty)$
be such that $C_{xy}=C_{yx}$. 
We call $C_{xy}$ the {\it conductance} between $x$ and $y$. 
We assume the following;

\medskip
\noindent (A1) There exist $c_1,c_2>0$ such that 
\[c_1\le \nu_x:=\sum_{y\in \Z^d}C_{xy}\le c_2\fal x\in \Z^d.\]

\medskip
\noindent (A2) There exist $M_0\ge 1,
\delta>0$ such that the following holds: 
for any $x,y\in \Z^d$ with $|x-y|=1$, there exist $N\ge 2$ and 
$z_1,\cdots, x_N\in B(x,M_0)$ such that $x_1=x$, $x_N=y$ 
and $C_{x_ix_{i+1}}\ge \delta$ for $i=1,\cdots, N-1$. 

\medskip
\noindent (A3) There exists  
a decreasing function $\vp: \mathbb{N}\to \R_+$ with $\sum_{i=1}^{\infty}i^{d+1}\vp(i)<\infty$
and $\vp(2i)\le c\;\vp (i)$ for all $i\in \mathbb{N}$ 
such that 
\[C_{xy}\le \vp (|x-y|)\qquad \mbox{for all }~~x, y\in \Z^d.\]

Note that (A1) and (A2) are very mild regularity conditions. 
(A1) prevents degeneracies, while (A2) says, roughly speaking, that the chain is locally irreducible in
a uniform way.
(A3) is the substantive assumption and says that the $C_{xy}$ satisfy a
uniform finite second moment condition. In fact, 
(A3) implies  the following:
\medskip
there exists $C_0>0$ 
such that 
\begin{equation}\label{eq:fin2mo}
\sup_{x\in \Z^d}\sum_{y\in \Z^d}|x-y|^2
C_{xy}\le C_0.
\end{equation} 
To see this, 
\begin{align}
\sum_{y\in \Z^d}|x-y|^2
C_{xy} & \le
\sum_{y\in \Z^d}|x-y|^2\vp(|x-y|)\label{eq:22a}\\
&=\sum_{i=0}^\infty\sum_{i< |x-y|\leq i+1} |x-y|^2 \vp(|x-y|)\nonumber \\
&\le c_3\sum_i (i+1)^2\vp (i) (i+1)^{d-1}
<\infty \nonumber
\end{align}
for all $x\in \Z^d$, where 
(A3) is used in the last inequality.

Define a symmetric Markov chain by
\[\P^x(X_1=y)=\frac {C_{xy}}{\nu_x}\fal 
x,y\in \Z^d.\]
Define $p_n(x,y):=\P^x(X_n=y)$ and $\bp_n(x,y)=p_n(x,y)/\nu_y$.
Note that $\bp_n(x,y)=\bp_n(y,x)$. By (A1), the ratio of $p_n(x,y)$
to $\bp_n(x,y)$ is bounded above and below by positive constants. 

Let $\mu_x\equiv 1$ for all $x\in \Z^d$ and for each $A\subset \Z^d$, define 
$\mu (A)=\sum_{y\in A}\mu_y=|A|$ and $\nu (A)=\sum_{y\in A}\nu_y$. 
Note that $L^2(\Z^d, \mu)=L^2(\Z^d, \nu)$
by (A1). Now, for each $f\in L^2(\Z^d,\mu)$, define 
\begin{eqnarray*}
\ce (f,f)&=&\tfrac 12 \sum_{x,y\in \Z^d}(f(x)-f(y))^2C_{xy},\\
\cf&=&\{f\in L^2(\Z^d,\mu): \ce(f,f)<\infty\}.
\end{eqnarray*}

It is easy to check $(\ce,\cf)$ is a regular Dirichlet form
on $L^2(\Z^d,\mu)$ and the generator is
\[\sum_{x,y\in \Z^d}(f(y)-f(x))
C_{xy}.\] 

Let $Y_t$ be the corresponding 
continuous time $\mu$-symmetric Markov chain on $\Z^d$. 
Let $\{U_i^x: i\in \N, x\in \Z^d\}$ be an independent sequence of exponential random 
variables, where the  parameter for $U_i^x$ is $\nu_x$, and  that is independent of $X_n$
and define $T_0=0, T_n=\sum_{k=1}^n U_k^{X_{k-1}}$. 
Define $T_0=0, T_n=\sum_{k=1}^n U_k$. 
Set $\widetilde Y_t=X_n$ if $T_n\le t<T_{n+1}$; 
it is well known that 
the laws of $\widetilde Y$ and $Y$ are the
same, and hence $\widetilde Y$ is a realization of the continuous time Markov
chain corresponding to (a time change of) $X_n$. 
Note that by (A1), the mean exponential holding time at each point 
for $\widetilde Y$ can be controlled uniformly from above and below by a positive 
constant. 
Let $p(t,x,y)$ be the transition density for $Y_t$ with 
respect to $\mu$.

We now introduce several processes related to $Y_t$, needed in what
follows. For each $D\ge 1$, let $\cs=D^{-1}\Z^d$ and define the 
rescaled process as $V_t=D^{-1}Y_{D^2t}$. Let $\mu^D$ be a measure
on $\cs$ defined by 
$\mu^D(A)=D^{-d}\mu (DA)=D^{-d}|A|$ for $A\subset \cs$. 
We can easily show that the Dirichlet form corresponding 
to  $V_t$ is 
\[\ce^D (f,f)=\tfrac 12 \sum_{x,y\in \cs}(f(x)-f(y))^2 D^{2-d}C_{Dx,Dy},\]
and the infinitesimal generator of $V_t$ is
\[ \sA^D f(x)=\sum_{y\in \cs} (f(y)-f(x) 
C_{Dx,Dy}D^2
=\sum_{x,y\in \cs}(f(y)-f(x))\frac{C_{Dx,Dy}D^{2-d}}{\mu^D_x},\]
for each $f\in L^2(\cs,\mu^D)$, where we denote $\mu^D_x:=\mu^D(\{x\})=D^{-d}$ for each $x\in \cs$. 
The heat kernel $p^D(t,x,y)$ for $V_t$ with respect to $\mu^D$ can be expressed as 
\begin{equation}\label{eq:heatrel}
p^D(t,x,y)=D^dp(D^2t, Dx,Dy)\fal x,y\in \cs, t>0.
\end{equation}
For $\lambda\ge 1$, let $W^{\lambda}_t$ be a process on $\cs$ 
with the  large jumps of $V_t$ removed. More precisely, $W^{\lambda}_t$ 
is a process whose Dirichlet form and infinitesimal
generator are 
\begin{eqnarray*}
\ce^{D,\lambda} (f,f)&=&\tfrac 12 \sum_{{x,y\in \cs}\atop
{|x-y|\le\lambda^{1/2}}}(f(x)-f(y))^2 D^{2-d}C_{Dx,Dy},\\
\calam f(x)&=&\sum_{{y\in \cs}\atop
{|x-y|\le\lambda^{1/2}}}(f(y)-f(x))\frac{C_{Dx,Dy}D^2}{\mu_{Dx}}.
\end{eqnarray*}
for each $f\in L^2(\cs,\mu^D)$. We denote the heat kernel for $W_t^\lambda$
by $p^{D,\lambda}(t,x,y)$, $x,y\in \cs$.

\section{Heat kernel estimates}
\subsection{Nash inequality}
For $f\in L^2(\Z^d, \mu)$, let
\[\ce_{NN}(f,f)=\tfrac 12 \sum_{{x,y\in \Z^d}; 
{|x-y|=1}}(f(x)-f(y))^2,\]
which is the Dirichlet form for the simple symmetric random walk in $\Z^d$.
We will prove the following Nash inequality.
\begin{propn} 
Assume (A2). 
There exists $c_1>0$ such that for any $f\in L^2(\Z^d, \mu)$,
\begin{equation}\label{eq:nash1}
\|f\|_2^{2(1+2/d)}\le c_1 \ce(f,f)\|f\|_1^{4/d}.
\end{equation} 
In particular, 
\begin{eqnarray}
p(t,x,y)\le c_1t^{-d/2}&\fal & x,y\in
\Z^d, t>0,\label{eq:nash2}\\
p^D(t,x,y)\le c_1t^{-d/2}&\fal & x,y\in\cs, t>0.\label{eq:nash3}
\end{eqnarray}
\end{propn}
\begin{rem}\label{remark1}
{\rm Since $p(t,x,y)=\P^x(Y_t=y)/\mu_y$, we have $p(t,x,y)\le 1/\mu_y$, 
so (\ref{eq:nash2}) is a crude estimate
for small $t$. However, we will continue to use it 
since we are mainly interested in the large time asymptotics.
}
\end{rem}
\medskip

\proof
Note that the equivalence of (\ref{eq:nash1}) and (\ref{eq:nash2})
is a  well-known fact (see \cite{CKS}).

The Markov chain corresponding to $\ce_{NN}$ is a 
(continuous time) simple random walk; let $r_t$ be its transition probabilities.
Since, as is well known, 
we have $r_t(x,x)\leq c t^{-d/2}$, then by \cite{CKS} we have 
\[\|f\|_2^{2(1+2/d)}\le c_1 \ce_{NN}(f,f)\|f\|_1^{4/d}
\fal f\in L^2(\Z^d, \mu).\]
See also \cite{SZ}.
By (A2), there exists $c_2>0$ such that  
\[\ce_{NN}(f,f)\le c_2\ce(f,f)\fal f\in L^2(\Z^d, \mu).\]
Using these facts and (\ref{eq:heatrel}), 
we have the desired result.\qed

\subsection{Exit time  probability estimates}

In this subsection, we will obtain some  exit time 
estimates. The argument presented here was first established in \cite{BL1} and
then extended and simplified in \cite{CK}, \cite{HK}. 

\begin{lem}\label{thm:kern2}
There exists $c_1>0$ such that
\begin{equation}\label{kern3}
p^{D,\lambda}(t,x,y)\leq c_1\;t^{-\frac{d}{2}}
\;\exp\left( -\lambda^{-\frac 12} |x-y|\right)
\end{equation}
for all $t\in(0,1]$,  
$x, y\in \cs$ and $\lambda\ge M_0$, where $M_0$ is given in (A2).
\end{lem}

\proof
Since $\lambda\ge M_0$, by (A2), we have 
$\ce_{NN}(f,f)\le c\ce^{D,\lambda} (f,f)$ for all $f\in L^2(\Z^d, \mu)$.
So we have (\ref{eq:nash1}) where 
$\ce(f,f)$ is replaced by $\ce^{1,\lambda} (f,f)$, and by a scaling argument 
we have 
\[p^{D,\lambda}(t,x,y)\le c_1t^{-d/2}\fal  x,y\in\cs, t>0.\]
Thus by Theorem (3.25) of \cite{CKS}, we have 
\begin{equation}\label{kern4}
p^{D,\lambda}(t,x,y)\leq c_1\;t^{-\frac{d}2}
\;\exp\left(-E(2t,x,y)\right)
\end{equation}
for all $t\le 1$ and $x,y\in \cs$, where  
\begin{eqnarray*}
E(t,x,y)&=&\sup\{|\psi(y)-\psi(x)|-t\; \Lambda(\psi)^2 :
\Lambda(\psi)<\infty\},\\
\Lambda(\psi)^2&=& \|e^{-2\psi}\Gamma_\lambda[e^\psi]\|_\infty \vee
\|e^{2\psi}\Gamma_\lambda[e^{-\psi}]\|_\infty,
\end{eqnarray*}
and $\Gamma_{\lambda}$ is defined by 
\begin{equation}  \label{densi}
\Gamma_{\lambda}[v](\xi)=
\sum_{{\eta,\xi\in\cs}\atop
{|\xi-\eta|\le\lambda^{1/2}}}(v(\eta)-v(\xi))^2
\frac{C_{D\eta,D\xi}D^2}{\mu_{D\xi}},\qquad 
\xi\in\cs.
\end{equation}
Now let $\psi(\xi)=\lambda^{-1/2}(|\xi-x|\wedge |x-y|)$.
Then, $|\psi(\eta)-\psi(\xi)|\le \lambda^{-1/2}|\eta-\xi|$,
so that 
\[(e^{\psi(\eta)-\psi(\xi)}-1)^2\le |\psi(\eta)-\psi(\xi)|^2
e^{2|\psi(\eta)-\psi(\xi)|}\le c\lambda^{-1}|\eta-\xi|^2\]
for $\eta,\xi\in \cs$ with $|\eta-\xi|\le \lambda^{1/2}$. 
Hence 
\begin{eqnarray*}
e^{-2\psi(\xi)}\Gamma_\lambda[e^\psi](\xi)&=&
\sum_{{\eta\in\cs}\atop
{|\xi-\eta|\le\lambda^{1/2}}}(e^{\psi(\eta)-\psi(\xi)}-1)^2
\frac{C_{D\eta,D\xi}D^2}{\mu_{D\xi}}\\
&\le &\lambda^{-1}\sum_{{\eta'\in \Z^d}
\atop{|\xi'-\eta'|\le D\lambda^{1/2}}}
|\eta'-\xi'|^2\frac{C_{\eta',\xi'}}{\mu_{\xi'}}\le C'
\end{eqnarray*}
for all $\xi\in \cs$ where (\ref{eq:fin2mo}) 
is used in the last inequality.
We have the same bound when $\psi$ is replaced by $-\psi$,
so $\Lambda(\psi)^2\leq {C'}^2$. 
Noting that $|\psi(y)-\psi(x)|\le \lambda^{-\frac 12}\; 
|x-y|$, we see that (\ref{kern3}) follows from (\ref{kern4}).
\qed

We now prove the following exit time estimate for the process.
For $A\subset \Z^d$ and a process $Z_t$ on $\Z^d$, let 
\[ \tau=\tau_A(Z):=\inf\{t\ge 0:
Z_t\notin A\}, \qquad T_A=T_A(Z):=\inf\{t\ge 0: Z_t\in A\}. \]
\begin{propn}\label{thm:4.2}
For $A> 0$ and $0<B<1$, there exists
$\gamma_i =\gamma_i(A,B)\in (0,1)$,
$i=1,2$, such that for every $D>0$ and $x\in
\Z^d$,
\begin{eqnarray}
\P^x \left( \tau_{B(x,\, AD)}(Y)<\gamma_1 \, D^2 \right)
&\le &B,\label{eq:tighty}\\
\P^x \left( \tau_{B(x,\, AD)}(X)<\gamma_2 \, D^2 \right)
&\le &B.\label{eq:tightx}
\end{eqnarray}
\end{propn}

\proof
It follows from  Lemma \ref{thm:kern2} that
for $t\in [1/4, \, 1]$ and $r >0$,
\begin{equation}\label{a6}
 \P^x \left( |W^{\lambda}_t-x|\geq r\right)
  = \sum_{ y\in \cs : \, |y-x|>r} p^{D,\lambda}(t, x, y) \mu^D_y  \leq c_1 \,I_{r,\lambda},
\end{equation}
where $I_{r,\lambda}:=e^{-\frac r2 \lambda^{-\frac 12}}$.
Define $\sigma_r:=\inf \{t\geq 0: \, 
|W_t^{\lambda}-W_0^{\lambda}|
\geq r \}$. Then by (\ref{a6}) and the strong Markov
property of $W^{\lambda}$ at time $\sigma_r$,
\begin{eqnarray*}
\P^x \left(\sigma_r \leq 1/2 \right)
&\leq & \P^x \left(\sigma_r \leq 1/2 \hbox{ and }
|W^{\lambda}_{1}-x|\leq r/2 \right)
 + \P^x  \left(|W^{\lambda}_{1}-x| > r /2 \right) \\
&\leq & \P^x \left(\sigma_r \leq 1/2 \hbox{ and }
|W^{\lambda}_{1}-W^{\lambda}_{\sigma_r}|
   ) > r /2 \right) + c_1\, I_{r/2,\lambda}  \\
&= & \mathbb{P}^x \left(1_{\{\sigma_r \leq 1/2\}}
\mathbb{P}^{W^{\lambda}_{\sigma_r}}\left( |W^{\lambda}_{1-{\sigma_r}}-
W^{\lambda}_{0}| > r /2 \right)\right) + c_1\, I_{r/2,\lambda}  \\
&\leq & \sup_{y\in B(x,r)^c}\sup_{s\le 1/2} \mathbb{P}^y
\left( |W^{\lambda}_{1-s}-y| > r /2 \right)
 +c_1\, I_{r/2,\lambda}\\
\end{eqnarray*}
Here in the second and the last inequalities, we used (\ref{a6}).
By the strong Markov property of $W^{\lambda}$,
for every $r > 0$,
\begin{eqnarray}
\P^x \left( \sup_{s\leq 1} |W^{\lambda}_{s}-
W^{\lambda}_0 | > r \right)
&\leq & \P^x (\sigma_r \leq 1/2) + \P^x (1/2 <\sigma_r \leq 1)
\nonumber\\
&\leq & c_2\, I_{r/2,\lambda} +
\mathbb{P}^x(\sigma_{r/2}\le 1/2)+
\mathbb{P}^x(\sigma_{r/2}> 1/2,\sigma_r \leq 1)\nonumber\\
&\leq & c_2\, I_{r/2,\lambda} +
\mathbb{P}^x(\sigma_{r/2}\le 1/2)+
\mathbb{E}^x \left[ \mathbb{P}^{W^{\lambda}_{1/2} } (\sigma_{r /2} \leq
1/2 ) \right]\nonumber\\
&\leq & c_3\, I_{r/4,\lambda}.\label{eq:kern5}
\end{eqnarray}
The constants $c_1, c_2, c_3>0$ above are independent of
$D \ge 1$, $x\in \cs$ and 
$\lambda \ge M_0$.

Now, define $B^{\lambda}$ to be the infinitesimal generator of $V_t$ with 
small jumps removed:
\begin{equation}\label{blambda}
B^{\lambda} v(\xi)=
\sum_{{\eta\in \cs}\atop
{|\eta-\xi|>\lambda^{1/2}}}(f(\eta)-f(\xi))\frac{C_{D\eta,D\xi}D^2}{\mu_{D\xi}}.
\end{equation}
Recall that $\calam$ is the generator of $W^\lambda$. 
We see that $\calam+B^\lambda$ is the generator for $V_t$.
Hence, if $Q^V_t$ and $Q^{W^\lambda}_t$
are the semigroups associated with $V_t$ and $W_t^\lambda$
respectively, we have that
\begin{equation}  \label{eq:semi}
Q^V_t v=Q^{W^\lambda}_tv +\sum_{k=1}^\infty S_k^{\lambda}(t)v,\qquad
v\in L^{\infty}(\cs,\mu^D),
\end{equation}
where
\begin{equation}\label{eq:sk}
S_k^{\lambda}(t)v=\int_0^t Q^{W^\lambda}_{t-s} B^{\lambda}
S_{k-1}^{\lambda}(s)v\, ds,\qquad k\geq1
\end{equation}
with $S_0^{\lambda}(t):=Q^{W^\lambda}_t$ (see, for example, 
Theorem $2.2$ in \cite{Le}).
Note that the series in (\ref{eq:semi}) defines a bounded linear 
operator on $L^\infty(\cs,\mu^D)$ for each $t>0$; this can be seen 
as follows. First, by (\ref{eq:22a}) and a simple calculation, we have
\begin{equation}\label{sekiin}
\sum_{{\eta\in \cs}\atop
{|\eta-\xi|>\lambda^{1/2}}}\frac{C_{D\eta,D\xi}D^2}{\mu_{D\xi}}
\le c_4\sum_{{y\in \Z^d}\atop
{|y-x|>D\lambda^{1/2}}}{C_{x,y}D^2}
\le \frac {c_5}{D^2\lambda}\sum_{y\in \Z^d}|x-y|^2\vp(|x-y|)D^2
\le \frac{c_6}{\lambda}.
\end{equation}
Using this, we see that there exists $c_7>0$
independent of $\lambda$ such that
$$
\|B^{\lambda} v\|_\infty\leq \frac{c_7}{\lambda}\|v\|_\infty.
$$
Noting that $\|Q^{W^\lambda}_tv\|_\infty\leq \|v\|_\infty$, 
by induction we have from (\ref{eq:sk}) that
\begin{equation}\label{sknorm}
\|S_k^{\lambda}(t)v\|_\infty\leq \frac{(c_8\lambda^{-1}
\;t)^k}{k!} \|v\|_\infty, \quad t>0,\; k\geq 1,
\end{equation}
and so the series above is bounded from
$L^\infty(\cs,\mu^D)$ to $L^\infty(\cs,\mu^D)$ for each $t>0$.

We will apply the above with $\lambda=M_0$.
By (\ref{sknorm}),
for any bounded function $f$ on $\cs$, we have
$$ \| Q^V_t f - Q^{W^\lambda}_t f \|_\infty
\leq \sum_{k=1}^\infty \frac{(c_8\lambda^{-1} \, t)^k} {k !} \, \| f \|_\infty
\leq c_9 \, t  \, e^{c_9 t } \, \| f \|_\infty.
$$
Applying this with $f$ equal to the indicator of $(\overline{B(\xi,r)})^c$, it
follows that there is a  constant $c_{10}>0$ that
is independent of $D \ge 1$ such that for every
$\xi\in \cs$ and every  $t\le 1$,
\begin{equation}\label{eqn:a11}
 \P^\xi \left( | V_t - \xi | > r \right) \leq 
 \P^\xi \left( | W^{M_0}_t -\xi | > r \right)
  + c_{10} \, t  .
\end{equation}
Applying the same argument we used in deriving
(\ref{eq:kern5}), we conclude there are positive constants
$c_{11}, c_{12}$ such that for $\xi\in \cs$,
\begin{equation}\label{eqn:a9}
\P^\xi \left( \sup_{s\leq t} | V_s - \xi|> r
\right) \leq c_{11} e^{-c_{12} r }+c_{11} \,t \qquad \hbox{for
every }  r >0 \hbox{ and  } t\leq 1.
\end{equation}
This implies that  for every $x\in \Z^d$,
$D'\ge 1$ and  $r >0$,
\begin{equation}\label{eqn:a10}
\P^x \left( \sup_{s\leq {D'}^2 \, t}| Y_s - x | >    r\, D'
\right) \leq c_{11} e^{-c_{12} r }+c_{11} \,t
\qquad \hbox{for every }  r >0 \hbox{ and  }
t\leq 1.
\end{equation}
For $A>0$ and $B\in (0, \, 1)$, we choose $r_0$ and $t_0$
so that $c_{11} e^{-c_{12} r_0 }+c_{11} \,t_0<B$
and take $D=r_0D'/A$. Then, by (\ref{eqn:a10}),
$$ \P^x \left( \sup_{s\leq  \gamma_1\,D^2 } |Y_s-x| \geq
A\,D \right) \leq B \qquad \hbox{for every }
D\ge r_0/A,
$$
where $\gamma_1=(A/r_0)^2t_0$. For $D< r_0/A$, we have
\begin{equation}\label{eq:r_0/a}
\P(U_1>\gamma_1\frac {r_0^2}{A^2})\le \P(U_1>\gamma_1D^2)\le 
\P^x \left( \sup_{s\leq  \gamma_1\,D^2 } |Y_s-x| <
A\,D \right),
\end{equation}
where $U_1$ is an exponential random variable with parameter 1. 
By (A1), the left hand side of (\ref{eq:r_0/a}) 
is greater than $1-B$ if $\gamma_1$ is taken to be small. 
Thus, (\ref{eq:tighty}) is proved. 

Now (\ref{eq:tightx}) can be proved in the same way as Theorem $2.8$ in 
\cite{BL1}.\qed

\section{Lower bounds and regularity for the heat kernel}

We now introduce the space-time process
$Z_s:=(U_s, V_s)$, where $U_s=U_0+s$. The filtration generated
by $Z$ satisfying the usual conditions
will be denoted by $\{ \widetilde \cf_s; \, s\geq 0\}$.
The law of the space-time
process $s\mapsto Z_s$ starting from $(t, x)$ will be denoted
as $\P^{(t, x)}$.
We say that  a non-negative Borel  measurable function
$q(t,x)$ on $[0, \infty)\times \cs$ is {\it parabolic}
in a relatively open subset $B$ of $[0, \infty)\times \cs$ if
for every relatively compact open subset $B_1$ of $B$,
$q(t, x)=\E^{(t,x)} \left[ q (Z_{\tau_{B_1}}) \right]$
for every $(t, x)\in B_1$,
where $\tau_{B_1}=\inf\{s> 0: \, Z_s\notin B_1\}$.

\medskip

We denote $\gamma:=\gamma(1/2, 1/2)<1$ the constant
in (\ref{eq:tighty}) corresponding to $A=B=1/2$.
For $t\ge 0$ and $r>0$, we define
\[  Q^D(t,x,r):=[t,t+\gamma r^{2}]\times (B(x,r)\cap \cs),
\]
where $B(x,r)=\{y\in \R^d: |x-y|\le r\}$. 

It is easy to see the following (see, for example, Lemma 4.5 in \cite{CK} for the proof).

\begin{lem}\label{4.6}
For each $t_0>0$ and $x_0\in \Z^d$, $q^D(t,x):=p^D (t_0-t, x,x_0)$
is parabolic on $[0,t_0)\times \cs$.
\end{lem}

The next proposition provides a lower bound for the heat kernel
and is the key step for the proof of the H\"older
continuity of $p^D (t, x, y)$. 

\begin{propn}\label{szhol} 
There exists $c_1>0$ and $\theta\in (0,1)$ such that if $|x-y|\leq t^{1/2}$,
$x,y\in \Z^d$ 
and $r>t^{1/2}/\theta$, then
\[ \P^x(Y_t=y, \tau_{B(x,r)}>t)\geq c_1t^{-d/2}. \]
\end{propn}

To prove this we first need some preliminary propositions.
A version of the following 
weighted Poincar\'e inequality can be found in  Lemma 1.19
of \cite{SZ}; we give an alternate proof.

\begin{lem}\label{1.19}
For $D\geq1$ and $l\in \Z^d$, let
\[ g_D(l)= c_0D^d \prod_{i=1}^d e^{-|l_i|/D},\]
where $c_0$ is determined by the equation $\sum_{l\in \Z^d}g_D(l)=1$. 
Then there exists $c_1>0$ such that 
\[c_1\Big<(f-\langle f \rangle_{g_D})^2\Big>_{g_D}\le D^{2-d}\sum_{l\in \Z^d}g_D(l)\sum_{i=1}^d(f(l+e^i)-f(l))^2,
\qquad  f\in L^2(\cs), \]
where 
\[\langle f \rangle_{g_D}=D^{-d}\sum_{l\in \Z^d}f(D^{-1}l)g_D(l)\]
and $e^i$ is the element of $\Z^d$ whose $j$-th component is $1$ if $j=i$ and $0$ otherwise.
\end{lem}

\def\dfm{{[f(m+1)-f(m)]}}
\def\dfn{{[f(n+1)-f(n)]}}

\proof A scaling argument shows that it suffices to consider only
the $D=1$ case.
Because of the product structure,
it is enough to consider the case when $d=1$. 

The weighted Poincar\'e inequality restricted to integers in $[-10,10]$,
i.e., where the sums are restricted to being over $\{-10, \ldots, 10\}$,
follows easily from the usual Poincar\'e inequality.
We will prove our weighted Poincar\'e inequality for positive
$k$ and the same argument works for negative $k$. These facts together
with the weighted Poincar\'e inequality on $[-10,10]$ and standard techniques as
in \cite{Je} give us the weighted Poincar\'e inequality for all of $\Z$.
So we restrict attention to nonnegative $k$.
Therefore all our sums below are over nonnegative integers.

Let
\begin{eqnarray*}
 I&=& \sum_{k,\ell} (f(k)-f(\ell))^2 e^{-k} e^{-\ell},\\
J_k&=& \sum_{\ell >k} \sum_{m=k}^{\ell-1} \sum_{n=k}^{\ell-1} \dfm\, \dfn e^{-\ell},\\
K&=&\sum_n \dfn^2 e^{-n}.\\
\end{eqnarray*}
Note
\[ I=\sum_k (f(k)-\langle{f}\rangle_{g_D})^2 e^{-k}, \]
so we need to show $I\leq c_2K$. 
We have, since $f(k)-f(\ell)=0$ when $k=\ell$,
\begin{eqnarray*}
I&=&2\sum_k \sum_{\ell>k} (f(k)-f(\ell))^2 e^{-k} e^{-\ell}\\
&=&2\sum_k \sum_{\ell>k} \Big(\sum_{m=k}^{\ell-1} \dfm\Big)\Big(\sum_{n=k}^{\ell-1} \dfn
\Big) e^{-k}e^{-\ell}\\
&=&2\sum_k J_k e^{-k}.
\end{eqnarray*}
We see that
\begin{eqnarray*}
J_k&=&\sum_{m\geq k}\sum_{n\geq k} \sum_{\ell> m\lor n} e^{-\ell} \dfm\, \dfn\\
&\leq &\sum_{m\geq k} \sum_{n\geq k} e^{-m\lor n} \dfm\, \dfn\\
&=&2\sum_{m\geq k} \sum_{n\geq m} e^{-n} \dfm\, \dfn \\
&=&2\sum_{n\geq k} \sum_{m=k}^{n-1} e^{-n} \dfm\, \dfn\\
&&\qquad\qquad+2 \sum_{n\geq k} e^{-n} \dfn^2\\
&\leq& 2\sum_{n\geq k} e^{-n} \dfn\, (f(n)-f(k))+2K.
\end{eqnarray*}
Hence
\begin{eqnarray*}
I&\leq& c_3\sum_k \sum_{n\geq k} e^{-n} \dfn\, (f(n)-f(k)) e^{-k} +\sum_k 2e^{-k}K\\
&\leq& c_4\Big(\sum_k \sum_{n\geq k} e^{-n} e^{-k} \dfn^2\Big)^{1/2}
\Big(\sum_k\sum_{n\geq k} e^{-n} [f(n)-f(k)]^2 e^{-k}\Big)^{1/2}\\
&&\qquad\qquad +c_4K\\
&\leq &c_4K^{1/2} I^{1/2} +c_4K.
\end{eqnarray*}
This implies 
\[ I\leq c_5K \]
as required.
\qed

The proof of the following lemma is similar to that of (1.16) in \cite{SZ}, but since
we need some modifications, we will give the proof.

\begin{lem}\label{1.20}
There is  an $\eps>0$ such that 
\begin{equation}\label{1.16}
p^D(t,D^{-1}k,D^{-1}m)\ge \eps t^{-d/2},
\end{equation}
for all $D\ge 1$, $(t,k,m)\in (D^{-1},\infty)\times \cs\times \cs$ with $|D^{-1}k-D^{-1}m|\le 2t^{1/2}$.
\end{lem}

\proof
First, note that it is enough to prove the following:
there is an $\eps>0$ such that 
\begin{equation}\label{1.25}
D^{-d}\sum_{l\in \Z^d}\log \Big(p^D(\tfrac 12, D^{-1}k,D^{-1}(l+m))\Big)g_D(l)\ge \tfrac 12 \log\eps,\end{equation} 
for all $D\ge 1$ and $k,m\in \Z^d$ with 
$|D^{-1}(k-m)|\le 2$. Indeed, by the Chapman-Kolmogorov equation, symmetry, and the fact 
$g_D(j)\le 1$ for all $j\in \Z^d$,
\[p^D(1,D^{-1}k,D^{-1}m)\ge D^{-d}\sum_jp^D(\tfrac 12, D^{-1}k,D^{-1}(j+k))p^D(\tfrac 12, D^{-1}m,
D^{-1}(j+k))g_D(j).\]
Thus, by Jensen's inequality, (\ref{1.25}) gives
\[p^D(1,D^{-1}k,D^{-1}l)\ge \eps\qquad D\ge 1, |D^{-1}k-D^{-1}l|\le 2.\] 
By a simple scaling argument, this gives (\ref{1.16}). 

So we will prove (\ref{1.25}). Set $u_t(l)=p^D(t,D^{-1}k,D^{-1}(l+m))$ and let
\[G(t)=D^{-d}\sum_{l\in \Z^d}\log(u_t(l))g_D(l).\]
By Jensen's inequality, we see that $G(t)\le 0$. Further,
\[G'(t)=D^{-d}\sum_{l\in \Z^d}\frac{\partial u}{\partial t}
(l)\frac{g_D(l)}{u_t(l)}
=-\ce^D(u_t(D\,\cdot),\frac{g_D(D\,\cdot)}{u_t(D\,\cdot)}).\]
Next, note that the following elementary inequality holds (see (1.23) of \cite{SZ} for the proof). 
\[\Big(\frac db -\frac ca\Big)(b-a)\le -\frac{c\wedge d}2 (\log b-\log a)^2+\frac{(d-c)^2}{2(c\wedge d)},\qquad
 a,b,c,d>0.\]
Hence
\begin{eqnarray*}
G'(t)&=&-\frac{D^{2-d}}2
\sum_{l\in \Z^d}\sum_{e\in \Z^d}\Big(\frac{g_D(l+e)}{u_t(l+e)}-\frac{g_D(l)}{u_t(l)}\Big)
\Big(u_t(l+e)-u_t(l)\Big)C_{l,l+e}\\
&\ge & \frac{D^{2-d}}2\sum_{l\in \Z^d}\sum_{e\in \Z^d}\frac{g_D(l+e)\wedge g_D(l)}2\Big(
\log u_t(l+e)-\log u_t(l)\Big)^2C_{l,l+e}\\
&&\ \ -\frac{D^{2-d}}2\sum_{l\in \Z^d}\sum_{e\in \Z^d}\frac{|g_D(l+e)-g_D(l)|^2}{2(g_D(l+e)\wedge g_D(l))}C_{l,l+e}\\
&\ge & cD^{2-d}\sum_{l\in \Z^d}\sum_{j=1}^d(g_D(l+e^j)\wedge g_D(l))
\Big(\log u_t(l+e^j)-\log u_t(l)\Big)^2\\
&&\ \ -D^{2-d}\sum_{l\in \Z^d}\sum_{e\in \Z^d}\frac{|g_D(l+e)-g_D(l)|^2}{4(g_D(l+e)\wedge g_D(l))}C_{l,l+e},
\end{eqnarray*}
where the last inequality is due to (A2) and the definition of $g_D$
(here recall that $e^i$ is in the element of $\Z^d$ whose $j$-th component is $1$ if $j=i$ and $0$ otherwise).
Note
$|g_D(l+e)-g_D(l)|\le c_1D^{-1}|e| (g_D(l+e)\wedge g_D(l))$.  Thus 
\begin{align*}
D^{2-d}&\sum_{l\in \Z^d}\sum_{e\in \Z^d}\frac{|g_D(l+e)-g_D(l)|^2}{4(g_D(l+e)\wedge g_D(l))}C_{l,l+e}
\le c_2D^{-d}\sum_{l}\sum_{e}C_{l,l+e}|e|^2 (g_D(l+e)\wedge g_D(l))\\
&\le  c_3\Big(\sup_{l}\sum_{e}C_{l,l+e}|e|^2\Big)\cdot D^{-d}\sum_{l}g_D(l)
= c_3\Big(\sup_{l}\sum_{e}C_{l,l+e}|e|^2\Big)<c_4,
\end{align*}
where we used (A3) in the last inequality. Note also $\min_{1\le i\le d} g_D(l+e^i)\ge c_5g_D(l)$.
Combining these, we have
\begin{eqnarray*}
G'(t)&\ge &c_6D^{2-d}\sum_{l\in \Z^d}\sum_{j=1}^d
\Big(\log u_t(l+e^j)-\log u_t(l)\Big)^2g_D(l)-c_4\\
&\ge & c_7D^{-d}\sum_l(\log u_t(l)-G(t))^2g_D(l)-c_4,
\end{eqnarray*}
where we used Lemma \ref{1.19} in the last inequality. 

Next, for $\sigma>0$, set $A_t(\sigma)=\{l\in \Z^d: u_t(l)\ge e^{-\sigma}\}$. Then, writing
$f^+$ and $f^-$ for the positive and negative parts of $f$, we have 
for each $\sigma>0$, 
\begin{align*}
D^{-d}&\sum_l(\log u_t(l)-G(t))^2g_D(l)\ge 
D^{-d}\sum_l(-(\log u_t)^-(l)-G(t))^2g_D(l)\\
&\ge  \frac{G(t)^2}{2D^d}\sum_{l\in A_t(\sigma)}g_D(l)-\sigma^2,
\end{align*}
where we used the elementary inequality $(A+B)^2\ge (A^2/2)-B^2$, $A,B\in \R$,  in the last inequality. 
Thus, we have
\begin{equation}\label{v2.v}
G'(t)\ge c_8I_{t,\sigma}G(t)^2-(c_4+\sigma^2),\end{equation}
where we let $I_{t,\sigma}=D^{-d}\sum_{l\in A_t(\sigma)}g_D(l)$. 
On the other hand, by (\ref{eq:tighty}) and scaling, we can find $r_0>2$ such that
\[D^{-d}\sum_{|D^{-1}l|\le r_0}p^D(t,D^{-1}k,D^{-1}(l+m))\ge 1/2,\qquad  D\ge 1, 
t\le 1, \mbox{ and } |D^{-1}(k-m)|\le 2.\]
In particular, if $\beta$ is the smallest value of  $e^{-2U}$ on $[-r_0,r_0]$, then
for each $t\in [1/4,1]$, 
\[1/2\le D^{-d}\sum_{|D^{-1}l|\le r_0}u_t(l)\le e^{-\sigma}r_0^d+(\sup_k |u_t(D^{-d}k)|)\cdot 
\frac {I_{t,\sigma}}{\beta^d}.\]
Thus by taking $\sigma=(4r_0^d)$ and using (\ref{eq:nash3}), we obtain 
$I_{t,\sigma}\ge c\beta^d$. Combining this with (\ref{v2.v}), there exists $0<\delta<1$ such that 
\begin{equation}\label{1.26s}
G'(t)\ge \delta G(t)^2-\delta^{-1},\qquad  D\ge 1, t\in [1/4,1], \mbox{ and }
|D^{-1}(k-m)|\le 2.\end{equation}

Now, by (\ref{1.26s}) and the mean value theorem, 
\begin{equation}\label{1.26ws}
G(1/2)-G(t)\ge -(4\delta)^{-1},\qquad  t\in [1/4,1].\end{equation}
We may assume $G(1/2)\le -5/(2\delta)$, since otherwise (\ref{1.25}) is clear. 
Then, by (\ref{1.26ws}) we have $G(t)\le -2\delta^{-1}$. So $\delta G(t)^2/2-\delta^{-1}\ge \delta^{-1}>0$.
So, by (\ref{1.26s}) again,
\[G'(t)\ge \delta G(t)^2/2,\qquad  t\in [1/4,1].\]
But this means that
\[G(\tfrac 12)^{-1}\le G(\tfrac 12)^{-1}-G(\tfrac 14)^{-1}=-\int_{1/4}^{1/2}\frac {G'(s)}{G^2(s)}ds
\le -\frac {\delta}8,\]
and therefore $G(1/2)\ge -8\delta^{-1}$. Thus (\ref{1.25}) holds with 
$\eps^{1/2}=\frac 12 \exp (-8\delta^{-1})$.\qed

\begin{lem}\label{offdiag}
Given $\delta>0$ there exists $\kappa$ such that if $x,y\in \Z^d$ and
$C\subset \Z^d$ with $\mbox{\rm dist}\,(x,C)$ and $\mbox{\rm dist}\,(y,C)$ both larger than
$\kappa t^{1/2}$, then 
\[ \P^x(Y_t=y, T_C\leq t)\leq \delta t^{-d/2}.
\]
\end{lem}

\proof By the strong Markov property we have
\begin{eqnarray*} \P^x(Y_t=y, T_C\leq t/2)&=& 
\P^x(1_{\{T_C\leq t/2\}}\P^{Y_{T_C}}(Y_{t-{T_C}}=y))\\
&\leq& c_1 (t/2)^{-d/2} \P^x(T_C\leq t/2).
\end{eqnarray*}
In Proposition \ref{thm:4.2} let us choose $A=1$ and $B=\delta/(4c_1 2^{d/2})$.
If we take $\kappa> (2\gamma_1)^{-1/2}$, then Proposition \ref{thm:4.2}
tells us that
\[ \P^x (T_C\leq t/2)\leq \P^x(\tau_{B(x,\kappa t^{1/2})}\leq t/2)\leq B,
\] and then
\begin{equation} \label{off1}
\P^x(Y_t=y, T_C\leq t/2)\leq \frac{\delta}{2} t^{-d/2}.
\end{equation}

We now consider
$ \P^x(Y_t=y, t/2\leq T_C\leq t).$
If the first hitting time of $C$ occurs between time $t/2$ and
time $t$, then the last hitting time of $C$ before time $t$ happens
after time $t/2$. So if $S_C=\sup\{s\leq t: Y_s\in C\}$, then
\[ \P^x(Y_t=y, t/2\leq T_C\leq t)\leq \P^x(Y_t=y, t/2\leq S_C\leq t). \]
We claim that by time reversal, 
\begin{equation}\label{offtime}  \P^x(Y_t=y, t/2\leq S_C\leq t) 
=\P^y(Y_t=x, T_C\leq t/2). 
\end{equation}
To see this, observe by the symmetry of the heat kernel $p$,
we have that if $t_i=(t/2)+ it/(2n)$, then
\begin{eqnarray*}
&&\P^x(Y_{t_k}=z_k,  \ldots, Y_{t_{n-1}}=z_{n-1}, Y_{t_n}=y)\\
&&=p(t_k,x,z_k) p(t/(2n),z_k,z_{k+1})\cdots p(t/(2n),z_{n-1},y)\\
&&=\P^y(Y_{t/(2n)}=z_{n-1}, \ldots, Y_{t-t_k}=z_k, Y_t=x).
\end{eqnarray*}
If we sum over $z_k\in C$ and $z_{k+1}, \ldots, z_{n-1}\notin C$, we have
\begin{align*}
\P^x( &Y_{t_k}\in C, Y_{t_{k+1}}\notin C, \ldots, Y_{t_{n-1}}\notin C, Y_t=y)\\
&= \P^y(Y_{t/(2n)}\notin C, \ldots,  Y_{t-t_{k+1}}\notin C, Y_{t-t_k}\in C, Y_t=x).
\end{align*}
If we sum over $k$, this yields
\[ \P^x(t/2\leq S_n'\leq t, Y_t=y)=\P^y(0\leq T_n'\leq t/2, Y_t=x), \]
where $S_n'=\sup\{t_k: Y_{t_k}\in C\}$ and $T_n'=\inf\{t_k: Y_{t_k}\in C\}$.
Letting $n\to \infty$ proves (\ref{offtime}).

Arguing as in the first part of the proof,
\[ \P^y(Y_t=x, T_C\leq t/2)\leq  \frac{\delta}{2} t^{-d/2}. \]
Therefore 
\[ \P^x(Y_t=y, t/2\leq T_C\leq t)\leq	 \frac{\delta}{2} t^{-d/2}, \] 
and combining with (\ref{off1}) proves the proposition.
\qed

\noindent 
{\sc Proof of Proposition \ref{szhol}.} 
We have from Lemma \ref{1.20} that there exists $\varepsilon$ such that
\[ p(t,x,y)\geq \varepsilon t^{-d/2} \]
if $|x-y|\leq 2t^{1/2}$. 
If we take $\delta=\varepsilon/2$ in Lemma \ref{offdiag}, then provided 
$r>\kappa t^{1/2}$, we have
\[ \P^x(Y_t=y, \tau_{B(x,r)}\leq t)\leq \frac{\eps}{2} t^{-d/2}. \]
Subtracting,
\[ \P^x(Y_t=y, \tau_{B(x,r)}>t)\geq \frac{\eps}{2} t^{-d/2} \]
if $|x-y|\leq 2t^{1/2}$, which is equivalent to what we want.
\qed

As a corollary of Proposition \ref{szhol} we have

\begin{cor}\label{szholcor}
For each $0<\eps<1$, there exists $\theta=\theta(\eps)>0$
with the following property:  if $D\geq 1$,
$x,y\in \cs$ with $|x-y|<t^{1/2}$, $r>0$, 
$t\in [0, (\theta r)^2)$, and  $\Gamma \subset B(y,t^{1/2})\cap \cs$ satisfies 
$\mu^D(\Gamma)t^{-d/2}\ge \eps$, then   
 \begin{equation}\label{ga55.2} 
\P^{ x} ( V_t\in \Gamma \mbox{ and } \tau_{B(x,r)}>t)\geq c_1\eps.\end{equation}
\end{cor}

\begin{lem}\label{krysaf} For each $0<\delta<1$, there exists 
$\gamma=\gamma_\delta\in (0,1)$ 
such that for $t>0$, $r>0$ and $x\in \cs$, if $A\subset 
Q^D_\gamma(t,x,r):=[t,t+\gamma_{\delta} r^{2}]\times (B(x,r)\cap\cs)$ 
satisfies $m\otimes \mu^D (A)/m\otimes \mu^D (Q^D_\gamma (t,x,r))\ge \delta$, then 
\[\P^{(t, x)} ( T_A(Z) < 
\tau_{Q^D_\gamma (t,x,r)}(Z)) \geq c_1\delta.\]
\end{lem}
\proof For each $\delta>0$, take $\gamma=\theta(\delta/4)^2$.  
Note that there exists $s=s_r\in [t+\delta\gamma r^2/4,t+\gamma r^2)$ such that 
\begin{equation}\label{elpar}
\mu^D (A_s)\ge \delta r^d/4\,\,\ge \frac {\delta}4 \Big(\frac{s-t}{\gamma}\Big)^{d/2}
\ge \frac {\delta}4 (s-t)^{d/2},
\end{equation}
where $A_s=\{(s,z)\in [0,\infty)\times \cs: (s,z)\in A\}$.  
Indeed, if not then 
\[m\otimes \mu^D (A)\le \delta\gamma r^{2+d}/4+(\gamma -\delta\gamma/4)\cdot (\delta/4)\cdot r^{2+d}
\le \delta\gamma  r^{2+d}/2,\]
which contradicts $m\otimes \mu^D (A)\ge \delta m\otimes \mu^D (Q^D_\gamma
(t,x,r))=\delta \gamma r^{2+d}$.
Now, using this fact and Corollary \ref{szholcor} (with $\eps=\delta/4$), we have
\begin{eqnarray*}
\P^{(t, x)} ( T_A(Z) < \tau_{Q^D_\gamma(t,x,r)}(Z)) &\geq & 
\P^{(t, x)} ( V_{s-t}\circ\theta_t\in A_s \mbox{ and } \tau_{B(x,r)}\circ\theta_t>s-t)\\
&\geq & c_1\delta/4, 
\end{eqnarray*}
which completes the proof.\qed
 
We will also use the following L\'evy system formula for $Y$ (cf. Lemma 4.7 in \cite{CK}).

\begin{lem}\label{4.9} 
Let $f$ be a non-negative measurable function on $\R_+\times \cs\times \cs$,
vanishing on the diagonal. Then for every $t\geq 0 $, $x\in \cs$ and 
a stopping time $T$ of $\{\cf_t \}_{t\geq 0}$,
\[  \E^x \left[\sum_{s\le T}f((s,V_{s-}, V_s)) \right]
=\E^x \left[\int_0^T \sum_{y\in \cs} f((s,V_s, y)) 
\frac{D^2C_{DV_s, Dy}}{\mu_{Y_{D^2s}}}\;ds \right]
\]
\end{lem}

Now we prove that the heat kernel $p^D (t, x, y)$ is H\"older
continuous in $(t, x, y)$, uniformly over $D$.
For $(t, x)\in [0, \infty)\times \cs$
and $r>0$ let $Q^D(t, x, r):=[t, \, t+\gamma r^2]
\times (B(x, R)\cap\cs)$,
where $\gamma:=\gamma(1/2, 1/2)\wedge \gamma_{1/3}<1$. 
Here $\gamma(1/2, 1/2)$ is the constant
in (\ref{eq:tighty}) corresponding to $A=B=1/2$ and 
$\gamma_{1/3}$ is the constant in Lemma \ref{krysaf} corresponding to
$\delta=1/3$.

The following theorem can be proved similarly to Theorem 4.1 in \cite{BL2} 
and Theorem 4.14 in \cite{CK}.
We will write down the proof for completeness.
\begin{thm}\label{4.15}
There are constants $c>0$ and $\beta>0$ (independent of $R,D$) such that for
every $0<R$, every $D\ge 1$, and every bounded parabolic function
$q$ in $Q^D(0, x_0, 4R)$,
\begin{equation}\label{eqn:holder1}
|q(s, x) -q(t, y)| \leq c \,  \| q \|_{\infty, R} \,
R^{-\beta} \, \left( |t-s|^{1/2} + |x-y|
\right)^\beta
\end{equation}
holds for $(s, x), \, (t, y)\in Q^D(0, x_0, R)$,
where $\| q \|_{\infty, R}:=\sup_{(t,y)\in [0, \, \gamma (4R)^2 ] \times \cs } |q(t,y)|$.
In particular, for the transition density function $p^D (t, x, y)$
of $V$,
\begin{equation}\label{eqn:holder2}
 |p^D (s, x_1, y_1) -p^D (t, x_2, y_2)| \leq c \,
t_0^{-(d+\beta)/2}
\left( |t-s|^{1/2} + |x_1-x_2|+|y_1-y_2| \right)^\beta,
\end{equation}
for any $0<t_0<1$, $t, \, s \in [t_0, \,  \infty)$ and $(x_i, y_i)\in
\cs\times \cs$ with $i=1, 2$. 
\end{thm}

\proof 
Recall that $Z_s=(U_s, V_s)$ is the space-time process of $V$,
where $U_s=U_0+s$. 
In the following, we suppress the
superscript $D$ from $Q^D(\cdot,\cdot,\cdot)$.
Without loss of generality, assume that
$0\leq q (z) \leq \| q \|_{\infty, R} =1$ for
$z\in  [0, \, \gamma\, (4R)^2 ] \times \cs$.
By Lemma \ref{krysaf}, there is a constant $c_1>0$
such that if $x\in \cs$, $0<r<1$ and $A\subset Q(t, x, r/2)$
with $\frac {m\otimes\mu^D (A)}{ m\otimes\mu^D (Q(t, x, r/2))} \geq 1/3$,
then
\begin{equation}\label{4.14}
 \P^{(t, x)} ( T_A(Z) < \tau_r(Z)) \geq c_1,
\end{equation}
where $\tau_r:=\tau_{Q(t, x, r)}$.
By Lemma \ref{4.9} with 
$f(s, y, z)=1_{B(x, r)}(y) \, 1_{\cs\setminus B(x, s)}(z)$ 
and $T=\tau_r$, 
there is a constant $c_2>0$ such that
if $s\geq 2 r$,
\begin{equation}\label{eqn:4.15}
\P^{(t, x)} (V_{\tau_r} \notin B(x, s) )
= \E^{(t, x)} \left[ \int_0^{\tau_r} 
\sum_{y\in\cs\setminus \overline{B(x, s)}}
\frac{D^2C_{DV_v, Dy}}{\mu_{Y_{D^2v}}}\, dv \right]
\leq \frac{c_2}{s^2} \E^{(t, x)}[\tau_r]\le \frac{c_2r^2}{s^2}.
\end{equation}
The first inequality of (\ref{eqn:4.15}) 
is due to the following computation.
\begin{eqnarray*}
\sup_{z\in B(x,r)\cap\cs}D^2\sum_{y\in\cs\setminus \overline{B(x, s)}}
C_{Dz,Dy}&\le &\sup_{z'\in B(Dx,Dr)}D^2\sum_{|z'-y'|\ge Ds/2}C_{z'y'}\le 
D^2\sum_{i>Ds/2}\vp (i)i^{d-1}\\
&\le &\frac 4{s^2}\sum_{i}\vp (i)i^{d+1}\le \frac c{s^2},\end{eqnarray*}
where (A3) is used in the last inequality. 
The last inequality of (\ref{eqn:4.15}) is due to the fact 
$\E^{(t, x)}[\tau_r]\le r^2$; this is clearly true since the time interval
for $Q(t,x,r)$ is $\gamma r^2$, which is less than $r^2$. 
($\E^x \tau_{B(x_0,r)}\leq c_1 r^2$ is also true -- see Lemma 
\ref{harnexp} (a).)
Let
$$ \eta= 1-\frac{c_1}4 \quad \hbox{and} \quad
\rho= \tfrac12 \wedge \left( \frac{\eta}2\right)^{1/2}
 \wedge \left( \frac{c_1 \, \eta} { 8 \, c_2} \right)^{1/2}.
$$
Note that for every $(t, x)\in Q(0, x_0, R)$, $q$ is parabolic in
$Q(t, x, R)\subset Q(0, \, x_0, \, 2R )$.
We will show that
\begin{equation}\label{eqn:4.16}
 \sup_{Q(t, x, \rho^k R )} q - \inf_{Q(t, x, \rho^k R )} q
 \leq  \eta^k \qquad \hbox{for all } k.
\end{equation}

For notational convenience, we write $Q_i$ for $Q(t, x, \rho^i R)$
and $\tau_i$ for $\tau_{Q(t, x, \rho^i R)}$. Define
$$ a_i = \inf_{Q_i} q \quad \hbox{and} \quad
  b_i =\sup_{Q_i} q.
$$
Clearly $b_i-a_i\leq 1 \leq \eta^i$ for all $i\leq 0$.
Now suppose that $b_i-a_i\leq \eta^i$ for all $i\leq k$
and we are going to show that $b_{k+1}-a_{k+1} \leq \eta^{k+1}$.
Observe that $Q_{k+1}\subset Q_k$ and so $a_k \leq q \leq b_k$ on $Q_{k+1}$.
Define
$$ A':=\{ z\in Q_{k+1}: \, q(z) \leq (a_k+b_k)/2 \}.
$$
We may suppose $\frac {m\otimes\mu^D (A')}{ m\otimes\mu^D (Q_{k+1} )} \geq 1/2$,
for if not we use $1-q$  instead of $q$.
Let $A$ be a compact subset of $A'$ such that
$\frac {m\otimes\mu^D (A)}{ m\otimes\mu^D (Q_{k+1} )} \geq 1/3$.
For any given $\varepsilon >0$, pick $z_1, z_2 \in Q_{k+1}$
so that $q(z_1)\geq b_{k+1}-\varepsilon$ and $q(z_2) \leq a_{k+1} + 
\varepsilon$. Then by (\ref{4.14})-(\ref{eqn:4.16}),
\begin{eqnarray*}
  b_{k+1} - a_{k+1}-2 \varepsilon 
&\leq&  q(z_1)- q(z_2)\\
&=&  \E^{z_1} \left[ q(Z_{T_A \wedge \tau_{k+1}})-q(z_2)
                   \right] \\
&=& \E^{z_1} \left[ q(Z_{T_A}) -q(z_2); \,
              T_A< \tau_{k+1}  \right]\\
&&\qquad + \E^{z_1} \big[ q(Z_{\tau_{k+1}})-q(z_2); \,
           T_A> \tau_{k+1}, \\
&&\qquad\qquad Z_{\tau_{k+1}}\in Q_k \big]\\
&&\qquad + \sum_{i=1}^\infty  \E^{z_1} \big[ q(Z_{\tau_{k+1}})-q(z_2); \,
           T_A> \tau_{k+1}, \\
&& \qquad\qquad  Z_{\tau_{k+1}  }\in Q_{k-i}
      \setminus  Q_{k+1-i} \big] \\
&\leq&  \left( \frac{a_k+b_k}2 -a_k \right) \P^{z_1} (  T_A< \tau_{k+1})\\
&&\qquad   + (b_k-a_k) \P^{z_1} ( T_A> \tau_{k+1} ) \\
&&\qquad  + \sum_{i=1}^\infty (b_{k-i}-a_{k-i}) \P^{z_1}
   ( Z_{\tau_{k+1}} \notin Q_{k+1-i}) \\
&\leq&  (b_k-a_k) \left( 1-\frac{\P^{z_1} (  T_A< \tau_{k+1})}2 \right)
  + \sum_{i=1}^\infty  c_2 \, \eta^k  (\rho^2
    / \eta )^i \\
&\leq&  ( 1-\frac{c_1}2)\, \eta^k  + 2c_2 \eta^{k-1} \rho^2  \\
&\leq&  ( 1-\frac{c_1}2) \eta^k + \frac{c_1}4 \eta^k \\
&=& \eta^{k+1}.
\end{eqnarray*}
Since
$\varepsilon$ is arbitrary, we have $b_{k+1}-a_{k+1}\leq \eta^{k+1}$
and  this
proves (\ref{eqn:4.16}).

For $z=(s, x)$ and $w=(t, y)$ in $Q(0, x_0, R)$ with $s\leq t$,
let $k$ be the largest integer
such that $|z-w|:= (\gamma^{-1} |t-s|)^{1/2}+|x-y| \leq \rho^k R$.
Then $\log (|z-w|/R) \geq (k+1) \log \rho $,
$w\in Q(s, x, \rho^k R)$ and
\begin{eqnarray*}
|q(z)-q(w)| \leq \eta^k = e^{k \log \eta} \leq
c_3 \left( \frac{|z-w|}{R} \right)^{\log \eta / \log \rho }.
\end{eqnarray*}
This proves (\ref{eqn:holder1}) with $\beta =
\log \eta / \log \rho$.

By (\ref{eq:nash2}) and Lemma \ref{4.6},
for every $0<t_0< 1$, $T_0\ge 2$ and $y \in \cs$,
$q(t, x):=p^D (T_0-t , x, y)$ is a parabolic function 
on $[0, T_0-\frac{t_0}2 ] \times \cs$ bounded above by 
$c_4 \, t_0^{-d/2}$.

For each fixed $t_0\in (0, \, 1)$ and $T_0\ge 2$, 
take $R$ such that $\gamma R^2 = t_0/2$.
Let $s, t\in [t_0, T_0]$ with $s>t$ and $x_1, x_2 \in \cs$.
Assume first that 
\begin{equation}\label{eqn:4.47}
|s-t|^{1/2}+|x_1-x_2| < 
\gamma^{1/2}\, R = (t_0/ 2)^{1/2}
\end{equation}
and so $(T_0-t, x_2)\in Q(T_0-s, x_1, R)\subset 
[0, T_0-\frac{t_0}2) \times \cs$. Applying
(\ref{eqn:holder1}) to the parabolic function
$q(t, x)$ with $(T_0-s, x_1)$, $(T_0-t, x_2)$ and
 $Q(T_0-s, x_1, R)$ in place of $(s, x)$, $(t, y)$
and $Q (0, x_0, R)$ there respectively, we have
\begin{equation}\label{eqn:holder3}
 |p^D (s , x_1, y) - p^D (t , x_2, y) | \leq c \, t_0^{-(d+\beta) /2}
   (|t-s|^{1/2} + |x_1-x_2| )^\beta.
\end{equation}
By (\ref{eq:nash3}), the inequality (\ref{eqn:holder3})
is true when (\ref{eqn:4.47}) does not hold.
So  (\ref{eqn:holder3}) holds for every
$t, s \in [t_0, T_0]$ and
$x_1, x_2 \in \cs$ for all $T_0\ge 2$. Inequality (\ref{eqn:holder2}) now
follows from (\ref{eqn:holder3}) by the symmetry of
$p(t, x, y)$  in $x$ and $y$. 
\qed

\section{Harnack inequality}

A function $h$ defined  on $\Z^d$ is harmonic on a subset $A$ of $\Z^d$ with respect to the Markov chain $X$ if
\[ \sum_z h(z) \P^x(X_1=z) =h(x), \qquad x\in A. \]
Because the Markov chain may not have bounded range, $h$ must be defined on all of $\Z^d$.
In order to avoid $h$ possibly being infinite in $A$, we will assume that $h$ is bounded on $\Z^d$,
but in what follows, the  constants do not depend at all on the $L^\infty$ bound on $h$. 
We say $h$ is harmonic with respect to $Y$ if $h(Y_{t\land \tau_A})$ is a $\P^x$-martingale
for each $x\in \Z^d$, where $\tau_A=\inf\{t: Y_t\notin A\}$.
It is not hard to see that a function is harmonic for $X$ if and only if it is harmonic for
$Y$, since the hitting probabilities of $X$
and $Y$ are the same. Also, because the  state space is discrete,
it is routine to see that a function is harmonic in a domain $A$ if and only if
$\ce(h,f)=0$ for all bounded $f$ supported in $A$; we will not use this latter fact.

In this section we first give an example of a symmetric random walk, i.e., where $\{X_{n+1}-X_n\}$ are
symmetric i.i.d.\ random variables, for which a uniform Harnack inequality fails.
Note that the Harnack inequality does hold for each ball of radius $n$,
but not with a constant independent of $n$.
Let $e^j$ be the unit vector in the $x_j$ direction, $j=1, \ldots, d$.

Let $b_n=n^{n^n}$ (or any other quickly growing sequence),
let $a_n$ be a sequence of positive numbers tending to 0, subject only
to $\sum a_n\leq 1/32$ and $\sum_n a_n b_n^2<\infty$. Let $\eps=\frac12 \sum a_n$.
Let ${\xi}_i$ be an i.i.d. sequence of random vectors
on $\Z^d$ with $\P^0({\xi}_1=\pm e^j)=(1-\eps)/(2d)$.
Let $\P^0({\xi}_1=\pm b_n e^1)=a_n$. Let $X_n=\sum_{i=1}^n {\xi}_i$.

Now let $\delta\in(0,1)$, $r_n=(1-\delta) b_n$, $z_n=(b_n,0)$,
$B_n=B(0,r_n)$, $\tau_n=\min\{k: X_k\notin B_n\}$, and $T_0=\min\{k:
X_k=0\}$. Define
$$h_n(x)=\P^x(X_{\tau_n}=z_n).$$
Each $h_n$ is a harmonic function in $B_n$.
If a uniform Harnack inequality were to hold, there would
exist $C$  not depending on $n$ such that
\[
h_n(0)/h_n(y)\leq C, \qquad y\in B(0, r_n/2).
\]

Since $\delta b_n \gg b_{n-1}$ for $n$ large, the only
way $X_{\tau_n}$ can equal $z_n$ is if the random walk jumps from
0 to $z_n$. So for $y_n\in B_n$, $y_n\ne 0$,
$$h_n(y_n)=\P^{y_n}(T_0<\tau_n)h(0).$$

But we claim that if $y_n\sim r_n/4$, then $\P^{y_n}(T_0<\tau_n)$ will
tend to 0  when $n\to \infty$, 
and then
$h_n(0)/h_n(y_n)\to \infty$. So no uniform Harnack inequality exists.

The claim is true is all dimensions greater than or equal to 2, but is easier
to prove when $d\geq 3$, so we concentrate on this case. We have
\begin{eqnarray*}
 \P^{y_n}(T_0<\tau_n) & \leq & 
 \P^{y_n}(T_0<\infty)=
\P^{y_n}(T_0<r_n^{1/4}) +\P^{y_n}(T_0\geq r_n^{1/4})\\
& \leq & \P^{y_n}(\max_{i\leq r_n^{1/4}} |X_i-X_0|\geq |y_n|)
+\sum_{i=[r_n^{1/4}]}^\infty \P^{y_n}(X_i=0).
\end{eqnarray*}
The first term on the last line goes to 0 by Doob's inequality (applied to each 
$(X_i, e^j)$, $j=1, \ldots, d$).
By Spitzer \cite{Spi}, p.~75, the sum above is bounded by
\[ c\sum_{i=[r_n^{1/4}]}^\infty \frac{1}{i^{d/2}}\leq c' (r_n^{1/4})^{1-(d/2)}, \]
which goes to 0 as $n\to \infty$.

Note that by taking $a_n$ tending to 0 fast enough, ${\xi}_1$ can
be made to be sub-Gaussian, or have even better tails.

Lawler \cite{Law} proved that the Harnack inequality
holds for a class of symmetric random walks with bounded range and also for a class of 
Markov chains with bounded range which are in general
not reversible. The content of the next proposition is that this continues
to be true for symmetric Markov chains with bounded range.

\begin{thm} \label{harnbdd}
Suppose the Markov chain has range bounded by $K$. 
Let $x_0\in \Z^d$. There exist constants
$c_1$ and $\theta$  not depending on  $x_0$ such that if 
$r\geq 4K(\theta^{-1}+1)$ and 
$h$ is nonnegative and bounded
on $\Z^d$ and harmonic on $B(x_0,r)$, then  
\[ h(x)\leq c_1 h(y), \qquad x,y\in B(x_0,\theta r). \]
\end{thm}

\proof First let us suppose that $d\geq 3$; we will remove this restriction at the
end of the proof. Let 
\[ G_B(x,y)=\E^x\int_0^{\tau_B} 1_{\{y\}}(Y_s)\, ds, \]
where $\tau_B=\inf\{t: Y_t\notin B(x_0,r)\}$.
$G_B$ is the Green function for the process $Y$ killed on exiting $B(x_0,r)$.
Since we are assuming $d\geq 3$ and $p(t,x,y)$ is always bounded by 
some constant, 
then by (\ref{eq:nash2}) we see that $G_B$ is bounded, say by $c_2$.

It follows by Proposition \ref{szhol} that there exists $\kappa$ such that $\P^x(Y_t=y, \tau_B>t)$ is bounded below by
$c_3t^{-d/2}$ provided $|x-x_0|, |y-x_0|\leq  t^{1/2}$ and $r>\kappa t^{1/2}$. 
Set $\theta=1/(4\kappa)$.
So integrating over $t\in [4\theta^2r^2,8\theta^2r^2]$, we see
$G_B(x,y)\geq c_5$ for $x,y\in B(x_0, 2\theta r)$.

Define $\overline h(x)=\E^x[h(Y_{T}); T<\tau_B],$ where $T=\inf\{t: Y_t\in B(x_0,2\theta r)\}$.
It is routine that $\overline h$ is equal to $h$ on $B'=B(x_0,\theta r)$, is 0 outside
of $B(x_0, r)$, and is excessive with respect to the process $Y_t$ killed
on exiting $B(x_0,r)$. 
The fact that $X$ has bounded range and $r>4K(\theta^{-1}+1)$ is what allows us to assert that 
$\overline h$ is equal to $h$ in $B(x_0, \theta r)$. 
See \cite{FOT}, p.~319, for the definition of excessive. 
By \cite{FOT}, Theorem 2.2.1,
there exists a measure $\pi$ supported on $\overline {B(x_0, r)}$
such that 
\[ \widetilde \ce (\overline h,v)=\int v(x) \, \pi(dx) \]
for all continuous $v$ with support contained in $B(x_0,r)$, 
where $\widetilde \ce$ is the Dirichlet form for $Y_t$ killed on exiting
$B(x_0,r)$. An easy approximation argumen shows that we also have
\[ \widetilde \ce(G_B \pi, v)=\int v(x) \, \pi(dx) \]
for such $v$, and we conclude $\overline h=G_B \pi$.
Since $\overline h$ is harmonic in $B(x_0, \theta r)$ and in 
$B(x_0,r)\setminus B(x_0, 2\theta r)$, it is not hard to see that $\pi$ in fact
is supported in $\overline{B(x_0,2\theta r)}\setminus B(x_0,\theta r)$. So for
$x,y\in B(x_0,\theta r)$, the upper and lower bounds on $G_B(x,y)$ imply
\begin{align*}
h(x)&=\sum_z G_B(x,z) \pi(\{z\}) \leq c_2 \pi(\Z^d)\\
&=\frac{c_2}{c_5} c_5\pi(\Z^d) \leq \frac{c_2}{c_5} \sum_z G_B(y,z) \pi(\{z\})\\
&=\frac{c_2}{c_5} h(y).
\end{align*}
This proves the theorem when $d\geq 3$.

When $d=2$, define a Markov chain $X'$ on $\Z^3$ by setting $C'_{(x_1,x_2,x_3),(y_1,y_2,y_3)}$
to be equal to $C_{(x_1,x_2), (y_1,y_2)}$ if $x_3=y_3$; 
equal to 1 if $x_1=y_1, x_2=y_2$, and $x_3= y_3 \pm 1$; and equal to 0 otherwise.
Suppose $h$ is harmonic with respect to $X$ on $A\subset\Z^2$. If we define $h'(x_1,x_2,x_3)=
h(x_1,x_2)$, it is routine to check that $h'$ is 
harmonic with respect to $X'$ on $A\times \Z\subset \Z^3$. The Harnack inequality
we just proved above applies to $h'$, and a Harnack inequality for $h$
then follows immediately.
\qed

As the example at the beginning of this section shows,
a uniform Harnack inequality need not hold when the
range is unbounded, so an additional assumption is needed to handle this case.
The assumption is modeled after \cite{BK}  and the proof is similar to the one in \cite{BL2}.
We assume
\medskip

(A4)\quad  There exists a constant $c_1$ such that $C_{xy}\leq c_1 C_{xy'}$ whenever
$|y-y'|\leq |x-y|/3$.

\begin{thm} \label{harnunbdd}
Suppose (A1)--(A3) hold and in addition (A4) holds. Suppose $x_0\in \Z^d$ and $R>M_0$, where
$M_0$ is defined in (A2).
There exists a constant $c_1$ such that if $h$ is nonnegative and
bounded on $\Z^d$ and harmonic on $B(x_0,2R)$, then
\begin{equation}\label{harin} 
h(x)\leq c_1 h(y), \qquad x,y\in B(x_0,R). 
\end{equation}
\end{thm}

Before proving Theorem \ref{harnunbdd} we prove a lemma. 
Note that (A4) is not needed for this lemma.

\begin{lem}\label{harnexp}
(a) $\E^x \tau_{B(x_0,r)}\leq c_1 r^2$.

(b) 
There exist $\theta\in (0,1)$ and $c_1, c_2>0$  such that
if $r> M_0/\theta$, then
$\P^x (\tau_{B(x_0,r)}\geq  r^2)\geq c_2$ 

and $\E^x \tau_{B(x_0,r)}\geq c_3 r^2$ if $x\in B(x_0,\theta r)$. 
\end{lem}  

\proof If $p(t,x,y)$ denotes the transition densities for $Y_t$, we know
\[ p(t,x,y)\leq c_4 t^{-d/2}. \]
So if we take $t=c_5r^2$ for large enough $c_5$, then
\[ \P^x(Y_t\in B(x_0,r))=\sum_{z\in B(x_0,r)} p(t,x,z)\leq c_6t^{-d/2} |B(x_0,r)| \leq \tfrac12. \]
This implies 
\[ \P^x(\tau_{B(x_0,r)}> t)\leq \tfrac12. \]
By the Markov property, for $m$ a positive integer
\[ \P^x(\tau_{B(x_0,r)}> (m+1)t)\leq \E^x[\P^{Y_{mt}}(\tau_{B(x_0,r)}> t); \tau_{B(x_0,r)}>mt]\leq \tfrac12 \P^x(\tau_{B(x_0,r)}>mt). \]
By induction,
\[ \P^x(\tau_{B(x_0,r)}>mt)\leq 2^{-m}, \]
and the first part of (a) follows. 

We also know
by Proposition \ref{szhol} that there exists $\kappa>1$ such that
\[ \P^x(Y_t=y, \tau_{B(x_0,r)}>t)\geq c_6 t^{-d/2} \]
if $|x-x_0|, |y-x_0|\leq t^{1/2}$ and $r>\kappa t^{1/2}$. 
Therefore taking $t=r^2/\kappa^2$,
\[ \P^x(\tau_{B(x_0,r)}>t)\geq \P^x(Y_t\in B(x_0,t^{1/2}), \tau_{B(x_0,r)}>t)\geq c_6 t^{-d/2} |B(x_0,t^{1/2})|\geq c_7 \]
if $x\in B(x_0, r/\kappa)$. Let $\theta=1/\kappa$. So $\E^x \tau_{B(x_0,r)}\geq t\P^x(\tau_{B(x_0,r)}>t)\geq c_7 r^2$, which proves (b).
\qed

\noindent {\sc Proof of Theorem \ref{harnunbdd}:}
Let $\kappa$ and $\theta $ be as in Lemma \ref{harnexp}.
 That a Harnack inequality inequality holds for each finite $R$ is
easy, provided $R\le 16M_0/\theta$,
so it suffices to assume $R>16 M_0/\theta$.

First of all, if $z_1\in \Z^d$ and  $w\notin B(z_1,2r)$, by
the L\'evy system formula,
\[ \E^x \sum_{s\leq \tau_{B(z_1,r)}\land t} 1_{(Y_{s-}\in B(z_1,r), Y_s=w)}=\E^x \int_0^{\tau_{B(z_1,r)}\land t} C_{Y_s,w}\, ds.
\]
Letting $t\to \infty$, we have
\[ \P^x(Y_{\tau_{B(z_1,r)}}=w)=\E^x\int_0^{\tau_{B(z_1,r)}} C_{Y_s,w}\, ds. \]
By (A4) the right hand side is bounded above by the quantity 
$c_{2} C_{z_1w}\E^x \tau_{B(z_1,r)}$ and below by  the quantity 
$c_{3} C_{z_1w}\E^x \tau_{B(z_1,r)}$.
By Lemma \ref{harnexp}, if $x,y\in B(z_1,\theta r)$, then $\E^x \tau_{B(z_1,r)}
\leq c_{4} \E^y \tau_{B(z_1,r)}$. We conclude
\[ \P^x(Y_{\tau_{B(z_1,r)}}=w)\leq c_{5} \P^y(Y_{\tau_{B(z_1,r)}}=w). \]
Taking linear combinations, if $H$ is a bounded function supported 
in $B(z_1,2r)^c$, then
\begin{equation} \label{3.2}
 \E^x H(Y_{\tau_{B(z_1, r)}})\leq c_{5} \E^y H(Y_{\tau_{B(z_1,r)}}), \qquad
x,y\in B(z_1, \theta r). 
\end{equation}

Choose $r_0=8M_0/\theta$. If $r\geq r_0$, then setting $t=r^2/\kappa^2$,
\[
\P^x(Y_t=y, \tau_{B(z_1,r)}>t)\geq c_{6} t^{-d/2}, \qquad x,y\in B(z_1, \theta r). 
\]
Summing over $A\subset B(z_1, \theta r)$, we see that
\begin{equation} \label{hhit}
 \P^x(T_A<\tau_{B(z_1,r)})\geq\P^x(Y_t\in A, \tau_{B(z_1,r)}>t)\geq c_{6} |A| t^{-d/2}
=c_{6}|A|r^{-d}, \qquad x\in B(z_1, \theta r).
\end{equation}
In particular, note that if $C\subset B(z_1,r)$ and
$|C|/|B(z_1,r)|\geq 1/3$, then
\begin{equation}\label{addlabel1}
\P^x(T_C<\tau_{B(z_1,r)})\geq c_{7}, \qquad x\in B(z_1, \theta r). 
\end{equation}

Next suppose $x,y\in B(z_1, \theta r_0)$. In view of (A2)
\[ \P^x(T_{\{y\}} <\tau_{B(z_1, r_0)})\geq c_{8}. \]
By optional stopping,
\begin{align*}
h(x)& \geq \E^x [h(Y_{T_{\{y\}}}); T_{\{y\}}<\tau_{B(z_1,r_0)}]\\
&=h(y)\P^x(T_{\{y\}}<\tau_{B(z_1,r_0)})\\
&\geq c_{8} h(y). 
\end{align*}

By looking at a constant multiple of $h$, we may assume
$\inf_{B(x_0,R)} h=1$. Choose
$z_0\in B(x_0,R)$ such that $h(z_0)= 1$. We want to show that $h$ is bounded above
in $B(x_0,R)$ by a constant not depending on 
$h$. 

Let
\begin{equation} \label{3.5}
\eta=\frac{c_{7}}{3}, \qquad \zeta=\frac13\land (c_{5}^{-1}\eta)\land c_{8}. 
\end{equation}

Now suppose there exists $x\in B(x_0,R)$ with $h(x)=K$ for some
$K$ large. Let $r$ be chosen so that
\begin{equation} \label{3.6}
2R^d/(c_{6}\zeta K)\leq 
|B(x_0,\theta r)|\leq 
4R^d/(c_{6}\zeta K).
\end{equation}
Note this implies
\begin{equation} \label{3.7}
r\leq c_9K^{-1/d}R.
\end{equation}
Without loss of generality we may assume $K$ is large enough that
$r\leq \theta R/4$.
 Let
\begin{equation} \label{3.75}
A=\{w\in B(x,\theta r): h(w)\geq \zeta K\}.
\end{equation}
By (\ref{hhit}) and optional stopping,
\begin{align*}
1\geq h(z_0)&\geq \E^{z_0}[h(Y_{T_A\land \tau_{B(x_0,2R)}});
T_A<\tau_{B(x_0,2R)}]\\
&\geq \zeta K \P^{z_0}(T_A<\tau_{B(x_0,2R)})\\
&\geq c_{6}\zeta K |A|R^{-d},\\
\end{align*}
hence
\begin{equation} \label{3.77}
\frac{|A|}{|B(x,\theta r)|}\leq 
\frac{R^d}{c_{6}\zeta K |B(x,\theta r)|}\leq \frac12.
\end{equation}
Let $C$ be a set contained in $B(x,\theta r)\setminus A$ such that
\begin{equation} \label{3.8}
\frac{|C|}{|B(x,\theta r)|}\geq \frac13. 
\end{equation}

Let $H=h1_{B(x,2r)^c}$.  We claim
$$\E^x[h(Y_{\tau_{B(x,r)}}); Y_{\tau_{B(x,r)}}\notin B(x,2r)]\leq \eta K.$$
If not
$$\E^x H(Y_{\tau_{B(x,r)}})> \eta K,$$
and by (\ref{3.2}), for all $y\in B(x,\theta r)$,
\begin{align*}
h(y)&\geq \E^y h(Y_{\tau_{B(x,r)}})\geq \E^y[h(Y_{\tau_{B(x,r)}}); Y_{\tau_{B(x,r)}}\notin B(x,2r)]\\
&\geq c_{5}^{-1} \E^x H(Y_{\tau_{B(x,r)}})\geq c_{5}^{-1}\eta K\\
&\geq \zeta K, \\   
\end{align*}
contradicting (\ref{3.8}) and the definition of $A$.

Let $N=\sup_{B(x,2r)} h(z)$. We then have
\begin{align*}
K&=h(x)=\E^x[h(Y_{T_{C}}); T_C<\tau_{B(x,r)}]+
\E^x[h(Y_{\tau_{B(x,r)}}); \tau_{B(x,r)}<T_{C}, Y_{\tau_{B(x,r)}}\in B(x,2r)]\\
&\qquad \qquad +\E^x[h(Y_{\tau_{B(x,r)}}); \tau_{B(x,r)}<T_{C}, Y_{\tau_{B(x,r)}}\notin B(x,2r)]\\
&\leq \zeta K\P^x(T_{C}<\tau_{B(x,r)})+N\P^x(\tau_{B(x,r)}<T_{C})+\eta K\\
&=\zeta K\P^x(T_C<\tau_{B(x,r)})+N(1-\P^x(T_C<\tau_{B(x,r)}))+\eta K,\\
\end{align*}
or
$$\frac{N}{K}\geq \frac{1-\eta-\zeta\P^x(T_C<\tau_{B(x,r)})}{1-\P^x(T_C<\tau_{B(x,r)})}.$$

Using (\ref{addlabel1})
there exists $\beta>0$ such that $N\geq K(1+\beta)$.
Therefore there exists $x'\in B(x,2r)$ with $h(x')\geq K(1+\beta)$.

Now suppose there exists $x_1\in B(x_0,R)$ with $h(x_1)=K_1$. Define
$r_1$ and $A_1$  in terms of $K_1$ analogously to (\ref{3.6})
and (\ref{3.75}). Using the above
argument (with $x_1$ replacing $x$ and $x_2$ replacing $x'$),
there exists $x_2\in B(x_1, 2r_1)$ with $h(x_2)=K_2\geq (1+\beta)K_1$.
We continue and obtain $r_2$ and $A_2$  and then $x_3,K_3,r_3, A_3$, etc. Note
$x_{i+1}\in B(x_i,2r_i)$ and $K_i\geq (1+\beta)^{i-1}K_1$. In view of
(\ref{3.7}), $\sum_i |x_{i+1}-x_i|\leq c_{10} K_1^{-1/d}R$. 
If $K_1$ is big enough, we  have a sequence $x_1,x_2, \ldots$ contained in $B(x_0,3R/2)$
Since $K_i\geq (1+\beta)^{i-1} K_1$ and $r_i\leq c_{11} K_i^{-1/d}R$,
there will be a first integer $i$ for which $r_i<2r_0$. But for all $y\in B(x_i, \theta r_i)$ we have
$h(y)\geq c_{8} h(x_i)$, so then $A_i= B(x_i, \theta r_i)$, a contradiction
to (\ref{3.77}).
\qed

\begin{cor} \label{jayrosen}
Let $\xi_i$ be an i.i.d.\ sequence of symmetric random vectors
taking values in $\Z^d$ with finite second moments. Let $X_n=\sum_{i=1}^n
\xi_i$ and suppose $X_n$ is aperiodic. 
Suppose there exists $c_1$ such that
\[ \P(\xi_1=y)\leq c_1 \P(\xi_1=y') \]
whenever $|y-y'|\leq |y|/3$. Then there exists $c_2$  and $R_0$ such
that for all $R$ larger than $R_0$ and any $w\notin B(x_0, R)$,
\[ \P^x(X_{\tau_{B(x_0,R)}}=w)\leq c_2 \P^y(X_{\tau_{B(x_0,R)}}=w),
\qquad x,y\in B(x_0,R/2). \]
\end{cor}

\proof We let $C_{xy}=\P(\xi_1=y-x)$. Since the $\xi_i$ are symmetric,
then the $X_n$ form a symmetric Markov chain, and it is easy
to see that (A1)--(A4) are satisfied.
We then apply Theorem \ref{harnunbdd} to $h(x)=\P^x(Y_{\tau_{B(x_0,R)}}=w)$.
\qed

\def\sE{{\cal E}}

\section{Central limit theorem}

\def\lam{{\lambda}}

Suppose we have a sequence $C_{xy}^n$ of conductances satisfying
(A1), (A2), and (A3) with constants and $\vp$ independent of $n$. Let $Y^{(n)}_t$ be the
corresponding continuous time Markov chains on $\Z^d$ and set
\[
Z^{(n)}_t=Y^{(n)}_{nt}/\sqrt n.
\] 

As noted previously, the Dirichlet form corresponding to the process $Z^{(n)}$
is
\begin{equation}\label{cltform}
\sE_n(f,f)=n^{2-d} \sum_{x,y\in n^{-1}\Z^d} (f(y)-f(x))^2 C^n_{nx,ny}.
\end{equation}
We will also need to discuss the form
\begin{equation}\label{cltform2}
\sE^R_n(f,f)=n^{2-d} \sum_{x,y\in n^{-1}\Z^d} (f(y)-f(x))^2 C^{n,R}_{nx,ny},
\end{equation}
where $C^{n,R}_{k,l}, k,l\in \Z^d$ is equal to $C^n_{k,l}$ if 
$|k-l|\leq nR$ and 0 otherwise.

Since the state space of $Z^{(n)}$ is $n^{-1}\Z^d$ while the limit process
will have $\R^d$ as its state space, we need to exercise some care with the
domains of the functions we deal with. First, if $g$ is defined on $\R^d$,
we define $R_n(g)$ to be the restriction of $g$ to $n^{-1}\Z^d$:
\[ R_n(g)(x)= g(x), \qquad x\in n^{-1}\Z^d. \]
If $g$ is defined on $n^{-1}\Z^d$, we next define an extension of $g$ to $\R^d$.
The one we use is defined as
follows. For $k\in \Z^d$, let 
\[ Q_n(k)=\prod_{j=1}^d [n^{-1}k_j, n^{-1}(k_j+1)]. \]
When $d=1$, we define the extension, $E_n(g)$, to be linear in each $Q_n(k)$ and to
agree with $g$ on the endpoints of each interval  $Q_n(k)$. For $d>1$ we define $E_n(g)$ inductively.
We use the definition in the $(d-1)$-dimensional case to define $E_n(g)$ on each face
of each $Q_n(k)$. We define $E_n(g)$ in the interior of a $Q_n(k)$ so that if 
$L$ is any 
line segment contained in the $Q_n(k)$ that is parallel to one of the
coordinate axes, then $E_n(g)$ is linear on $L$. For example, when $d=2$, $n=1$, and $k=(0,0)$, then
\begin{align*}
E_n(g)(s,t)&=g(0,0)(1-s)(1-t)+g(0,1)(1-s)t+g(1,0)s(1-t)\\
&\qquad +g(1,1)st, \qquad \qquad 0\leq s,t\leq1. 
\end{align*}

Recall $e^j$ is the unit vector in the $x_j$ direction and let $(x,y)$ denote the inner
product in $\R^d$. 
If $k=(k_1, \ldots, k_d)\in \Z^d$, let $\sP(k)$ be the union
of the line segment from 0 to $(k_1, 0, \ldots, 0)$, the line
segment from $(k_1, 0, \ldots, 0)$ to $(k_1, k_2, 0, \ldots, 0)$,
..., and the line segment from $(k_1, \ldots, k_{d-1},0)$ to
$k$. 
For $z\in \Z^d$ and $1\leq i\leq d$, let
\[
L^i_z=\{(y,k)\in (n^{-1}\Z^d)^2: y+ n^{-1}\sP(nk) \mbox{ contains the line segment from $z$ to }
z+n^{-1}e^i\}.
\]
We note that $(x,k)\in L^i_z$ for $z\in n^{-1}\Z^d$ if and only if 
$(x+k)_l=z_l$ for $l=1,...,i-1$, $x_l=z_l$ for $l=i+1,...,d$ and $z_i\in [x_i\wedge 
(x+k)_i, x_i\vee (x+k)_i)$. So, for each $k$, the number of $x$ that satisfies
$(x,k)\in L^i_z$ is at most $n|k_i|$.

Recall $\sgn r$ is equal to 1 if $r>0$, equal to 0 if 
$r=0$, and equal to $-1$ if $r<0$. 
We define a map $a^n$ from $\R^d$
into $\cal M$, the collection  of $d\times d$ matrices as follows:
Fix $R$.  If $x\in n^{-1}\Z^d$, let
the $(i,j)$-th element of $a^n$  be given by
\begin{equation}\label{defan}
\big(a^n(x)\big)_{ij}
=\sum_{(y,k)\in L^i_x}   C_{ny, n(y+k)}^{n,R} nk_{j}\,\sgn k_i.
\end{equation}
For general $x=(x_i)_{i=1}^d\in\R^d$, we define $a^n(x):=a^n([x]_n)$, where 
we set $[x]_n=(n^{-1}[nx_i])_{i=1}^d$. 
$a^n(x)$ is not symmetric in general, but under (A5), we see that 
$(a^n(x))_{ij}$ is bounded for all $i,j$, (which can be proved similarly
to (\ref{aolacop}) below) and when $n$ is large, we can use Cauchy-Schwarz, etc.,
as in  the symmetric case. 
Note that if $C_{xy}^n=0$ for $|x-y|>1$ (i.e., the nearest neighbor case), then 
the expression in (\ref{defan}) is equal to 
$2C_{nx,nx+e^i}^n$ if $i=j$ and equal to 0 if $i\ne j$.
(In particular, $a^n(x)$ is symmetric in this case.)

We make the following assumption.

\medskip
\noindent (A5) There exist $R>0$ and a Borel measurable $a:\R^d\to\cal M$ such that 
$a$ is symmetric and uniformly elliptic, 
the map $x\to a(x)$ is continuous, and 
$a^n$ converges to $a$ uniformly on compacts sets. 
\medskip

We will see from the proofs below that if (A5) holds for one $R$, then it holds
for every $R>1$ and the limit $a$ is independent of $R$.

Since 
$a$ is uniformly elliptic, if we
define
\[\ce_a(f,f)=\int_{\R^d}(\nabla f(x), a(x)\nabla f(x))dx,\]
then $(\ce_a, H^1(\R^d))$ is a regular Dirichlet form on $L^2(\R^d, dx)$ where
$H^1(\R^d)$ is the Sobolev space of square integrable functions with one square integrable derivative.
Further, it is well-known that the corresponding heat kernel $p^a(t,x,y)$ satisfies the following
estimate,
\begin{equation}\label{aroest}
c_1t^{-d/2}\exp\Big(-c_2\frac{|x-y|^2}t\Big)\le p^a(t,x,y)\le
c_3t^{-d/2}\exp\Big(-c_4\frac{|x-y|^2}t\Big),
\end{equation}
for all $t>0$ and all $x,y\in \R^d$. 
As a consequence, the corresponding diffusion (which we denote by $\{Z_t\}$)
can be defined without ambiguity from any starting point.

In this section we prove the following central limit theorem.
Let $C([0,t_0]; \R^d)$ be  the collection of continuous paths from $[0,t_0]$ to $\R^d$.

\begin{thm}\label{clt}
Suppose (A1)-(A3) and (A5) hold. 

(a) Then for each $x$ and each $t_0$ the  
$\P^{[x]_n}$-law of $\{Z^{(n)}_t; 0\leq t\leq t_0\}$ converges weakly with respect to the topology of
the space $D([0,t_0], \R^d)$. The limit probability gives full  measure to $C([0,t_0], \R^d)$.

(b) If $Z_t$ is the canonical process on $C([0,\infty), \R^d)$ and $\P^x$ is the weak limit
of the $\P^{[x]_n}$-laws of $Z^{(n)}$, then
the process $\{Z_t, \P^x\}$ has continuous paths and is the symmetric process corresponding 
to the Dirichlet form $\ce_a$. 
\end{thm}

Before giving the proof, we discuss three examples. First, suppose each $X^{(n)}$ is the
sum of i.i.d.\ random vectors. Then the $C^{n}_{xy}$ will depend only
on $y-x$, and so the $a^n(x)$ will be constant in the variable $x$. Therefore, if convergence holds, the limit
$a(x)$ will be constant in $x$. 
This means that the limit is a linear transformation of $d$-dimensional Brownian motion, as one
would expect.

For another example, suppose the $X^{(n)}$ are nearest neighbor Markov chains, i.e., $C^{n}_{xy}=0$
if $|x-y|\ne 1$. Then in this case the result of \cite{SZ} is included 
in our Corollary \ref{d1corsz} and \ref{cltcor}.

Third,  suppose $C^n_{xy}=C_{xy}$ does not depend on $n$. Unless $C_{xy}$ is a function only
of $y-x$, then (2.6) of \cite{SZ} (which is (\ref{szstrassmp}) below) 
will not be satisfied, and this situation is covered by
Theorem \ref{clt} but not by the results of \cite{SZ}. To be fair, the goal of \cite{SZ}
was not to obtain a general central limit theorem, but instead to come up with a way of
approximating diffusions by Markov chains.
Condition (A5) is restrictive. 
For this $C^n_{xy}=C_{xy}$ case, if we further assume that $C_{xy}=0$ for $|x-y|>1$,
then $a(x)$ is always a constant matrix. Indeed, in this case
the expression in (\ref{defan}) is equal to 
$2C_{nx,nx+e^i}\delta_{ij}$, which converges
to $(a(x))_{ij}$ uniformly on compacts
as $n\to\infty$ by (A5). So, for any $m\in \N$, the limit of $a^n(x/m)$ is equal to
$a(x)$, i.e., $a(x/m)=a(x)$. Since $a$ is continuous, we conclude $a(x)=a(0)$ for all $x\in \R^d$. 

Before we prove Theorem \ref{clt}, we prove a proposition showing tightness of the laws of $Z^{(n)}$.

\begin{propn}\label{clttight}
Suppose $\{n_j\}$ is a subsequence. Then there exists a further subsequence $\{n_{j_k}\}$
such that

(a) For each $f$ that is $C^\infty$ on $\R^d$  with compact support, $E_{n_{j_k}}(P_t^{n_{j_k}} R_{n_{j_k}}(f))$ converges
uniformly on compact subsets;  if we denote the limit by $P_t f$, then the operator
$P_t$ is linear and 
extends to all continuous functions on $\R^d$ with compact support and is the semigroup of a symmetric
strong Markov process on $\R^d$ with continuous paths.

(b) For each $x$ and each $t_0$ the $\P^{[x]_{n_{j_k}}}$ law of $\{Z_t^{(n_{j_k})}; 0\leq t \leq t_0\}$ converges weakly to a probability
$\P^x$ giving full measure to $C([0,t_0]; \R^d)$.
\end{propn}

\proof Let $t_0>0$ and $\eta>0$. Let $\tau_n$ be stopping times 
bounded by $t_0$ and let $\delta_n\to 0$. Then by 
Proposition \ref{thm:4.2} and the strong Markov property, 
\[
\limsup_{n \to \infty }\P( |{ Z}^{(n)}_{\tau_n+\delta_n}-{ Z}^{(n)}_{\tau_n}|>\eta)=0.
\]
This, Proposition \ref{thm:4.2}, and \cite{A} imply that
the laws of the $\{{Z}^{(n)}\}$ are tight in $D[0,t_0]$ for each $t_0$. 

Fix $t_0$ and $\eta>0$. ${Z}^{(n)}$ will have a jump of size larger than $\eta$
before time $t_0$ only if $|Y^{(n)}_{t}-Y^{(n)}_{t-}|\geq \eta\sqrt n$
for some $t\leq nt_0$. 
By the L\'evy system formula, the probability of this is bounded by
\begin{align*}
\E^x\sum_{s\leq nt_0} 1_{(|Y^{(n)}_s-Y^{(n)}_{s-}|\geq \eta\sqrt n)}&=\E^x \int_0^{nt_0} \sum_{|x-Y^{(n)}_s|\geq
\eta\sqrt n} C^n_{Y^{(n)}_sx}\, ds\\
&\leq c_1(nt_0) \sum_{i\geq \eta \sqrt n} \vp(i) i^{d-1}\\
&\leq c_1 t_0 \eta^{-2} \sum_{i\geq \eta \sqrt n} \vp(i) i^{d+1}, 
\end{align*}
which tends to 0 by dominated convergence as $n\to \infty$. Since
this is true for each $t_0$ and $ \eta>0$
we conclude that any subsequential limit point of the sequence ${Z}^{(n)}$
will have  continuous paths.

{}From this point on the argument is fairly standard. We give a sketch,
leaving the details to the reader. Take a 
countable dense subset $\{t_i\}$ of $[0,\infty)$ and a countable dense
subset $\{f_m\}$ of the $C^\infty$ functions on $\R^d$ with compact support.
Let $P_t^n$ be the semigroup for $Z^{(n)}$. In view of Theorem \ref{4.15},
$E_{n_{j}}(P_{t_i}^{{n_{j}}}(R_{n_{j}}(f_m)))$ will be equicontinuous. By a diagonalization argument,
we can find a subsequence $\{n_{j_k}\}$ of $\{n_j\}$ such that for 
each $i$ and $m$, as $n_{j_k}\to \infty$,
these functions converge uniformly on compact
sets. Call the limit $P_{t_i} f_m$. Using the equicontinuity, we can 
define $P_t f_m$ by continuity for all $t$, and because the norm of each $P_t$ is bounded by 1, we
can also define $P_t f$ by continuity for $f$ continuous with compact support.
Using the equicontinuity yet again, it is easy to see that the $P_t$ satisfy the 
semigroup property and that $P_t$ maps continuous functions with compact support
into continuous functions. One can thus construct a strong Markov process that
has $P_t$ as its semigroup. The symmetry of $P^{(n)}_t$ leads to the symmetry of 
$P_t$.  

For each $x$, the $\P^{[x]_{n_j}}$ laws of $\{Z^{(n_{j})}_t; 0\leq t\leq t_0\}$
are tight. Fix $x$, let $\{n'\}$ be any subsequence of $\{n_{j_k}\}$
along which the $\P^{[x]_{n'}}$ converge weakly,
and let $\P$ be the weak limit of the subsequence
$\P^{[x]_{n'}}$. Suppose $F$ is a continuous
functional on $C([0,t_0]; \R^d)$ of the form $F(\omega)=\prod_{\ell=1}^L g_i(\omega(s_i))$, 
where the $g_i$ are continuous with compact support and $0\leq s_1<\cdots < s_L\leq t_0$.
When $L=1$, then
\begin{align*}
\E g_1(Z_{s_1})&=\lim \E^{[x]_{n'}} R_{n'}(g_1)(Z_{s_1}^{(n')})\\
&=\lim P_{s_1}^{n'} R_{n'}(g_1)([x]_{n'})\\
&=P_{s_1} g_1(x).
\end{align*}
Thus the one-dimensional  distributions of a  subsequential limit point of the $\P^{[x]_{n_{j_k}}}$
do not depend on the subsequence $\{n'\}$.
Using the Markov property of $Z^{(n)}$ and the equicontinuity, a similar argument
shows that the same is true of
 the $L$-dimensional distributions.
 Therefore there must be weak convergence along the subsequence
$\{n_{j_k}\}$.
As proved above, the weak limit is concentrated on the set of continuous paths.
\qed

\noindent{\sc Proof of Theorem \ref{clt}:} We denote the Dirichlet form for the
process $Z^{(n)}$ by $\ce_n$. 
Suppose $f, g$ are $C^\infty$ on $\R^d$ with compact support.
Let $U^n_\lam$ be the $\lam$-resolvent for $Z^{(n)}$; this means that
\[ U_\lam^n h(x)=\E^x\int_0^\infty e^{-\lam t} h(Z^{(n)}_t)\, dt \]
for $x\in n^{-1}\Z^d$ and $h$ having domain $n^{-1}\Z^d$. We write
$P_t^n$ for the semigroup for $Z^{(n)}$. 

Using Proposition \ref{clttight}, we need to show that if we have a subsequential
limit point of the $P^n_t$ in the sense of that proposition, then the limiting
process corresponds to the Dirichlet form $\ce_a$. 
Let $\{n'\}$ be a subsequence of $\{n\}$
for which the subsequence converges in the sense of Proposition \ref{clttight}, and let $U_\lam$ be the 
$\lam$-resolvent of the limiting process. 

Let $F_{n'}=U^{n'}_\lam (R_{n'}(f))$. Then
\begin{equation}\label{cltA} 
\ce_{n'}(F_{n'}, R_{n'}(g))=(R_{n'}(f), R_{n'}(g))-\lam (F_{n'}, R_{n'}(g)), 
\end{equation}
where we let 
$(h_1,h_2)=\sum_{x\in n^{-1}\Z^d} h_1(x)h_2(x)\mu^D_x$ 
for functions defined on $n^{-1}\Z^d$. (Recall that our base measure is $\mu^D$.)  
Let $H_n=E_n(F_n)$ and $H=U_\lam f$. 
The equicontinuity result of Theorem \ref{4.15} and Proposition \ref{clttight}
shows that the $H_{n'}$ converges uniformly on compacts to $H$.
If we can show
\begin{equation}\label{clt1}
\ce_a(H,g)=(f,g)-\lam (H,g),
\end{equation}
this will show that the $\lam$-resolvent for the limiting process is the
same as the $\lam$-resolvent for the process corresponding to $\ce_a$, and
the proof will be complete; we also use $(h_1,h_2)$ to denote $\int h_1(x) h_2(x)\, dx$
when $h_1, h_2$ are functions defined on $\R^d$.

Next, since $f\in L^2(\R^d)$ and $f$ is $C^\infty$, then $R_n(f)\in L^2(d\mu_n)$.
Standard Dirichlet form theory shows that
\[ \| U^n_\lam (R_n(f))\|_2\leq \frac{1}{\lam}\| R_n(f)\|_2, \]
that is, the $L^2$ norm of $F_n$ is bounded in $n$.
We see that   
\begin{equation}\label{clt2.5}
 \int |\nabla H_n(x)|^2 \, dx \leq c_1 \ce_n(F_n,F_n)=c_1((R_n(f),F_n)-\lam(F_n,F_n)) 
\end{equation}
is bounded in $n$. By the compact imbedding of $W^{1,2}$ into $L^2$, we conclude
that $\{H_n\}$ is a compact sequence in $L^2(\R^d)$; here
$W^{1,2}$ is the space of functions whose gradient is square integrable. Since $H_{n'}$ converges
on compacts to $H$, it follows that $H_{n'}$ converges in $L^2$ to $H$. 
We note also that by (\ref{cltA})
\begin{equation}\label{clt2.7}
\ce_n(F_n,F_n)=(R_n(f),F_n)-\lambda (F_n, F_n)
\end{equation}
is uniformly bounded in $n$.

We need to know that 
\begin{equation}\label{clt2}
|\sE^{R}_n(F_n,R_n(g))-\ce_{a^n}(H_n,g)|\to 0
\end{equation}
as $n\to \infty$.
The proof of this is a bit lengthy and we defer it to Lemma \ref{cltlem} below.

We also need to show that 
\begin{equation}\label{clt2q}
|\sE_n(F_n,R_n(g))-\sE^{R}_n(F_n,R_n(g))|\to 0
\end{equation}
as $n\to \infty$. This follows because by Cauchy-Schwarz,  
we have
\begin{align*}
\Bigl| \sum_{x,y\in n^{-1}\Z^d}& (F_n(y)-F_n(x))
n^{2-d}C^{n}_{nx,ny} 
(R_n(g)(y)-R_n(g)(x))\\
&-\sum_{x,y\in n^{-1}\Z^d} (F_n(y)-F_n(x))
n^{2-d} C^{n,R}_{nx,ny} (R_n(g)(y)-R_n(g)(x))\Bigr|\\
\leq & c\ce_n(F_n,F_n)^{1/2}
\Big[\sum_{x,y\in n^{-1}\Z^d} 
n^{2-d}(C^{n}_{nx,ny}-C^{n,R}_{nx,ny}) (R_n(g)(y)-R_n(g)(x))^2\Big]^{1/2}.
\end{align*}
The term within the brackets on the last line  is bounded by
\[
 c\| \nabla g\|^2_\infty 
\sup_{x\in n^{-1}\Z^d}\sum_{y\in \Z^d, |x-y|>nR} |x-y|^2 C^n_{xy}
\leq c'\sum_{i>nR} i^{d-1} i^2 \vp(i),
\]
which will be less than $\eps^2$ if $n$ is large. 

Using (\ref{cltA}), (\ref{clt1}), (\ref{clt2}), and (\ref{clt2q}), we see that it suffices to
show
\begin{equation}\label{clt3}
\ce_{a^{n'}}(H_{n'},g)\to \ce_a(H,g). 
\end{equation}
Now 
\begin{equation}\label{clt4} 
|\ce_{a^{n'}}(H_{n'},g)-\ce_a(H_{n'},g)|=\Bigl|\int \nabla H_{n'} \cdot (a^{n'}-a) \nabla g\Bigr|. 
\end{equation}
Since $\nabla g$ is bounded with compact support and $|\nabla H_{n'}|$ is bounded in $L^2$, then
(A5) and the Cauchy-Schwarz inequality tell us 
that the right hand side of (\ref{clt4}) tends to 0 as $n\to \infty$.
Therefore we need to show
\begin{equation}\label{clt5}
\ce_a(H_{n'},g)\to \ce_a(H,g).
\end{equation}
But if $\nabla h$ is bounded with compact support,
then
\begin{equation}\label{clt6}
 \int (\nabla H_{n'})\, h=-\int H_{n'} \nabla h\to -\int H\, \nabla h=\int (\nabla H)\, h. 
\end{equation}
If we take the supremum over such $h$ that also have $L^2$ norm bounded
by 1, then Fatou's lemma and the Cauchy-Schwarz inequality show
that $\nabla H$ is in $L^2$.
If $h$ is bounded with compact support, let $\varepsilon>0$ and
approximate $h$  by a $C^1$ function $\tilde h$ 
with compact support  such that $\|h-\tilde h\|_2\leq \varepsilon$.
Since  $|\nabla H_n|$ is bounded in $L^2$, 
then $|\int \nabla H_{n'}(h-\tilde h)|\leq c_1\varepsilon$
and $|\int \nabla H(h-\tilde h)|\leq c_1\varepsilon$.
So by (\ref{clt6})
\[ \limsup_{n'\to \infty} \Bigl|\int \nabla H_{n'}h-\int \nabla H \, h\bigr|\leq 2c_1\varepsilon. \]
Because $\varepsilon $ is arbitrary, we have
\begin{equation}\label{clt7}
 \int \nabla H_{n'}\, h\to \int \nabla H\, h. 
\end{equation}
If we apply (\ref{clt7}) with $h=a\nabla g$, we obtain (\ref{clt5}).
\qed

To complete the proof we have

\begin{lem}\label{cltlem}
With the notation of the above proof,
\[ |\sE^{R}_n(F_n,R_n(g))-\ce_{a^n}(H_n,g)|\to 0
\]               
as $n\to \infty$.
\end{lem}
 
\proof
{\it Step 1.} Let $\eps, \eta_1, \eta_2, \delta>0$ and let $\{\sS_m\}$ be a collection of 
cubes with disjoint interiors whose union
contains the support of $g$ and such that the oscillation of $a$ on each $\sS_m$ is less than
$\eta_1$ and the oscillation of $\nabla g$ on each $\sS_m$ is less than $\eta_2$. One way to
construct such a collection is to take a cube large enough to contain the support of $g$,
divide it into $2^d$ equal subcubes, and then divide each of the subcubes and so on
until the oscillation restrictions are satisfied.

{\it Step 2.}
Let $\sS_m'$ be the cube with the same center as  $\sS_m$ but side length $(1-2\delta)$
times as long. Let $A=\cup_m(\sS_m-\sS_m')$. We claim it suffices to show that
\begin{align}
\Bigl| \int_{A^c} &\nabla H_n(x) \cdot a^n(x) \nabla g(x)\, dx\nonumber\\
&- \sum_{x\notin A, x\in n^{-1}\Z^d}\ \sum_{y\in n^{-1}\Z^d} (F_n(y)-F_n(x)) 
n^{2-d}C^{n,R}_{nx,ny} (R_n(g)(y)-R_n(g)(x))\Bigr|\nonumber\\
&\qquad\to 0 \label{t7A}
\end{align}
as $n\to \infty$.
To see this, note first that
by Cauchy-Schwarz and (\ref{clt2.5})
\begin{align*}
 \int_{A} \nabla H_n(x) \cdot a^n(x) \nabla g(x)\, dx
&\leq \ce_{a^n}(H_n,H_n)^{1/2} \Big(\int_{A} \nabla g(x) \cdot a^n(x) \nabla g(x)\, dx\Big)^{1/2}\\
&\leq c\ce_{a^n}(H_n,H_n)^{1/2} \| \nabla g\|_\infty  |A|^{1/2}
\end{align*} will be less than $\eps$  if $\delta$ is taken sufficiently small. 
Next note that for any $x\in n^{-1}\Z^d$,
\begin{align*} \sum_{y\in n^{-1}\Z^d}n^{2-d}C^{n,R}_{nx,ny} (R_n(g)(y)-R_n(g)(x))^2
&\leq n^{-d} \|\nabla g\|_\infty^2  
\sum_{y\in n^{-1}\Z^d} C^{n,R}_{nx,ny} |ny-nx|^2\\
&\leq cn^{-d}.
\end{align*}
So by Cauchy-Schwarz and (\ref{clt2.7})
\begin{align} \sum_{x\in A, x\in n^{-1}\Z^d} \sum_{y\in n^{-1}\Z^d} &(F_n(y)-F_n(x)) 
n^{2-d}C^{n,R}_{nx,ny} (R_n(g)(y)-R_n(g)(x))\nonumber\\
&\leq \ce_n (F_n,F_n)^{1/2} \Big(
\sum_{x\in A, x\in n^{-1}\Z^d}\sum_{y\in n^{-1}\Z^d}n^{2-d}C^{n,R}_{nx,ny} 
(R_n(g)(y)-R_n(g)(x))^2\Big)^{1/2}\nonumber\\
&\leq  c\ce_n (F_n,F_n)^{1/2} \Big(n^{-d}\, \mbox{\rm card}\, (A\cap n^{-1}\Z^d)\Big)^{1/2},
\label{crdel}
\end{align}
which will be less than $\eps$ if $\delta $ is taken small enough and $n$ is large.

{\it Step 3.} Let $x_m$ be the center of $\sS_m$. Define $\ol g$ by requiring $\ol g$
to be linear on each $\sS_m$ and satisfying $\ol g(x_m)=g(x_m)$, $\nabla\ol g(x_m)=\nabla g(x_m)$.
We claim it suffices to show that
\begin{align}
\Bigl| \int_{A^c}& \nabla H_n(x) \cdot a^n(x) \nabla \ol g(x)\, dx\nonumber\\
&- \sum_{x\notin A, x\in n^{-1}\Z^d}\ \sum_{y\in n^{-1}\Z^d} (F_n(y)-F_n(x)) 
n^{2-d}C^{n,R}_{nx,ny} (R_n(\ol g)(y)-R_n(\ol g)(x))\Bigr|\nonumber\\
&\qquad \to 0 \label{t7B}
\end{align}
To see this, note that
\begin{align*}
\Bigl| \int_{A^c} \nabla H_n(x) \cdot a^n(x) \nabla \ol g(x)\, dx-
 &\int_{A^c} \nabla H_n(x) \cdot a^n(x) \nabla  g(x)\, dx \Bigr|\\
&\leq \ce_{a^n} (H_n, H_n)^{1/2} 
\Big( \int_{A^c} \nabla (\ol g-g)(x) \cdot a^n(x) \nabla (\ol g-g)(x)\, dx\Big)^{1/2}\\
&\leq c\ce_{a^n} (H_n, H_n)^{1/2} \eta_2,
\end{align*}
which will be less than $\eps$ if $\eta_2$ is chosen small enough.
A similar argument shows that the difference between the second term in (\ref{t7B}) and the corresponding
term with $\ol g$ replaced by $g$ is small; cf. Step 2.

{\it Step 4.} Let $\ol C^n_{xy}=C^{n,\delta/2}_{xy}$ and define
$\ol a^n(x)$ by $(\ol a^n(x))_{ij}
=\sum_{(y,k)\in L^i_x}  \ol C_{ny, n(y+k)}^{n} nk_{j}\sgn k_i$.
We claim it suffices to show that 
\begin{align}
\Bigl| \int_{A^c}& \nabla H_n(x) \cdot \ol a^n(x) \nabla \ol g(x)\, dx\nonumber\\
&- \sum_{x\notin A, x\in n^{-1}\Z^d}\ \sum_{y\in n^{-1}\Z^d} (F_n(y)-F_n(x))
n^{2-d}\ol C^{n}_{nx,ny} (R_n(\ol g)(y)-R_n(\ol g)(x))\Bigr|\nonumber\\
&\qquad \to 0 \label{t7C}
\end{align}
To prove this, we first note that the following can be proved in the same way as 
(\ref{clt2q}).

\begin{align*}
\Bigl|\sum_{x\notin A, x\in n^{-1}\Z^d}& \sum_{y\in n^{-1}\Z^d} (F_n(y)-F_n(x))
n^{2-d}\ol C^{n}_{nx,ny} (R_n(\ol g)(y)-R_n(\ol g)(x))\\
&-\sum_{x\notin A, x\in n^{-1}\Z^d}\ \sum_{y\in n^{-1}\Z^d} (F_n(y)-F_n(x))
n^{2-d} C^{n,R}_{nx,ny} (R_n(\ol g)(y)-R_n(\ol g)(x))\Bigr|\to 0,
\end{align*}
as $n\to\infty$. Next,
\begin{align}
\Bigl| \int_{A^c}& \nabla H_n(x) \cdot \ol a^n(x) \nabla \ol g(x)\, dx
-\int_{A^c} \nabla H_n(x) \cdot a^n(x) \nabla \ol g(x)\, dx\Bigr|\nonumber\\
&\le c\Bigl(\int_{A^c} (\nabla H_n(x))^2dx\Bigr)^{1/2}
\Bigl(\int_{A^c}(\ol a^n(x)- a^n(x))(\nabla \ol g(x))^2\Bigr)^{1/2}.\label{aola}
\end{align}  
We can estimate
\begin{align}
\Bigl|(\ol a^n(x)- a^n(x))_{ij}\Bigr|&\le \sum_{(y,k)\in L^i_x}   
 |\ol C_{ny, n(y+k)}^{n}-C_{ny, n(y+k)}^{n}|nk_{j}\nonumber\\
&\le  c_1\sup_{x\in \Z^d}\sum_{y\in \Z^d, |x-y|>n\delta/2} |x-y|^2 C^n_{xy}
\leq c_2\sum_{i>n\delta/2} i^{d-1} i^2 \vp(i),\label{aolacop}
\end{align}
where in the second inequality, 
we used the fact that for each $k$, the number of 
$y$ that satisfies $(y,k)\in L^i_z$ 
is at most $n|k_i|$ (as mentioned when we defined $L^i_z$).
So the right hand side of (\ref{aola}) will be less than 
$\eps$ if $n$ is large. 

{\it Step 5.} We have chosen the $\sS_m$ so that the oscillation of $a$ on each $\sS_m$
is at most $\eta_1$. Since we have that the $a^n$ converge to the $a$ uniformly
on compacts and there are only finitely many $\sS_m$'s, then for $n$ large the oscillation
of $a^n$ on any $\sS_m$ will be at most $2\eta_1$. 

{\it Step 6.} We will now prove (\ref{t7C}). By Step 3, $\ol g$ is linear on each 
$\sS_m$, so it is enough to discuss the case where $\ol g(x)=x_{j_0}$ on $\sS'_m$ 
for some $j_0$ and then use a linearity argument. 
Noting that $H_n=F_n$ on $n^{-1}\Z^d$, define
\[\ce_n^{\sS'_m}(H_n,\ol g):=
\sum_{x\in \sS'_m\cap n^{-1}\Z^d}\ 
\sum_{y\in n^{-1}\Z^d} (H_n(y)-H_n(x))
n^{2-d}\ol C^{n}_{nx,ny} (R_n(\ol g)(y)-R_n(\ol g)(x)).\]
Since there is no term involving different $\sS'_m$'s, we will consider each $\sS'_m$
separately.
We will fix an $x_0\in \sS'_m$ and look at the terms
involving $H_n(x_0+n^{-1} e_i)-H_n(x_0)$. First, by an elementary computation using 
the definition of the linear extension map $E_n$, we have 
\begin{equation}\label{t71}
\int_{Q_n(x_0)}\frac{\partial H_n}{\partial x_i}dx=\frac 1{2^{d-1}n^{d-1}}\sum_{z\in
V_i(x_0)}(H_n(z+n^{-1}e_i)-H_n(z))
\end{equation}
where $V_i(x_0)$ is the collection of vertices of the face of $Q_n(x_0)$ perpendicular to
$e_i$ and with the smaller $e_i$ component. (E.g., for a square, $V_1(x_0)$ is the two
leftmost corners, $V_2(x_0)$ is the two lower corners.) So
\begin{align*}
\int_{Q_n(x_0)}(\nabla H_n,\ol a^n\nabla \ol g)dx&=\sum_{i,j=1}^d\int_{Q_n(x_0)}
\frac{\partial}{\partial x_i} H_n \ol a^n_{ij}\frac{\partial}{\partial x_j} \ol g\, dx=
\sum_i\ol a^n_{ij_0}(x_0)\int_{Q_n(x_0)}\frac{\partial}{\partial x_i} H_n dx\\
&=\sum_{i=1}^d\ol a^n_{ij_0}(x_0)\frac 1{2^{d-1}n^{d-1}}\sum_{z\in
V_i(x_0)}(H_n(z+n^{-1}e_i)-H_n(z)).
\end{align*}
Summing over all cubes that contains $H_n(x_0+n^{-1}e^i)-H_n(x_0)$, the coefficient in front of 
$H_n(x_0+n^{-1}e^i)-H_n(x_0)$ will be 
\begin{equation}\label{t72}
\frac{n^{1-d}}{2^{d-1}}
\sum_{z\in V_i(x_0+n^{-1}e^i-e_*)}\ol a^n_{ij_0}(z),
\end{equation}         
where $e_*=(1/n,...,1/n)$.

We next look at $\ce_n^{\sS'_m}(H_n,\ol g)$.  
Since $\ol g(x+k)-\ol g(x)=k_{j_0}$ where $k=(k_1,...,k_d)$, we have 
\[
\ce_n^{\sS'_m}(H_n,\ol g)=n^{2-d}
\sum_{{x\in \sS'_m\cap n^{-1}Z^d,}\atop{k\in n^{-1}Z^d}} 
(H_n(x+k)-H_n(x)) 
\ol C^{n}_{nx,n(x+k)}k_{j_0}.
\]
Let us fix $x$ and $k$ and replace $(H_n(x+k)-H_n(x))$ by the sum $\sum_{m=1}^{|k|} 
(H_n(z_{m+1})-H_n(z_m))$ (here $|k|:=|k_1|+...+|k_d|$ and $|z_{m+1}-z_m|=1/n$)
so that the union of the line segments belongs to $x+ n^{-1} \sP(k)$. 
We will get a term of the form $H_n(x_0+n^{-1}e_i)-H_n(x_0)$
if $z_m=x_0$ and $z_{m+1}=x_0+n^{-1}e_i$ (we get $H_n(x_0)-H_n(x_0+n^{-1}e_i)$ if 
$z_{m+1}=x_0$ and $z_{m}=x_0+n^{-1}e_i$), so the contribution will be
\[ n^{2-d}\ol C^n_{nx,n(x+k)}k_{j_0}(\sgn~ k_i). \]
Summing over $x\in \sS'_m\cap n^{-1}Z^d,k\in n^{-1}\Z^d$, we have that the coefficient in front of 
$H_n(x_0+n^{-1}e^i)-H_n(x_0)$ for $\ce_n^{\sS'_m}(H_n,\ol g)$ is 
\begin{equation}\label{t73}
n^{2-d}\sum_{{x\in \sS'_m\cap n^{-1}Z^d,}\atop {(x,k)\in L^i_{x_0}}} 
\ol C^n_{nx,n(x+k)}k_{j_0}(\sgn~ k_i).
\end{equation}
On the other hand, by the definition of 
$\ol a^n$, we have 
\begin{equation}\label{t73q}
n^{2-d}\sum_{(x,k)\in L^i_{x_0}} \ol C^n_{nx,n(x+k)}k_{j_0}(\sgn~ k_i)=n^{1-d}
\ol a^n_{ij_0}(x_0).
\end{equation}
Let $\sS_m''$ be the cube with the same center as  $\sS_m'$ but side length $(1-2\delta)$
times as long. If $x_0\in \sS_m''\cap n^{-1}\Z^d$, then the expressions in (\ref{t73}) and (\ref{t73q}) are equal, because
$\ol C^n_{nx,n(x+k)}=0$ for 
$x\notin \sS_m'\cap n^{-1}\Z^d, (x,k)\in L^i_{x_0} $.
Since the oscillation of $a^n$ on each $\sS'_m$ is less that $2\eta_1$ as in Step 5, 
by (\ref{aolacop}) the oscillation of $\ol a^n$ on each $\sS'_m$ is less that $3\eta_1$.
Thus, when $x_0\in \sS_m''\cap n^{-1}\Z^d$, we see
that the absolute value of the difference between (\ref{t72}) and (\ref{t73}) is bounded by  $3\eta_1n^{1-d}$. 
(Note that $\mbox{\rm card}\, V_i(x_0+n^{-1}e^i-e_*)=2^{d-1}$ is used here.)
Now, if $x_0\in (\sS'_m-\sS_m'')\cap n^{-1}\Z^d$, then
the difference between (\ref{t72}) and (\ref{t73}) is bounded by
$c_*n^{1-d}$, because similarly to (\ref{aolacop}) we have 
\[\sum_{(x,k)\in L^i_{x_0}} \ol C^n_{nx,n(x+k)}nk_{j_0}(\sgn~ k_i)
\le c_1\sup_{x\in \Z^d}\sum_{y\in \Z^d} |x-y|^2 C^n_{xy}
\leq c_2\sum_{i} i^{d-1} i^2 \vp(i)=:c_*.\]
Denote $H_{x_0,i}:=H_n(x_0+n^{-1}e^i)-H_n(x_0)$, $A':=(\cup_m(\sS_m'-\sS_m''))\cap n^{-1}\Z^d$ and 
$B:=(\cup_m\sS_m'')$$\cap n^{-1}\Z^d$.
Using the Cauchy-Schwarz inequality, we have 
\begin{eqnarray}
&\Big|\int_{\cup_m\sS'_m}(\nabla H_n,\ol a^n\nabla \ol g)dx-\sum_m\ce_n^{\sS'_m}(H_n,\ol g)\Big|
\label{compcdis}\\
\le &\eta_1n^{1-d}\sum_{x_0\in B,i=1,\cdots, d}|H_{x_0,i}|
+c_*n^{1-d}\sum_{x_0\in A',i=1,\cdots, d}|H_{x_0,i}|\nonumber\\
\le& c_1\eta_1\Big(n^{-d}\mbox{\rm card}\, B\Big)^{1/2}
\Big(n^{2-d}\sum_{x_0\in n^{-1}\Z^d,i}(H_{x_0,i})^2\Big)^{1/2}
\\&~~~~+c_*\Big(n^{-d}\mbox{\rm card}\, A'\Big)^{1/2}
\Big(n^{2-d}\sum_{x_0\in n^{-1}\Z^d,i}(H_{x_0,i})^2\Big)^{1/2}\nonumber\\
\le &c_2(\eta_1+\eps)\Big(n^{2-d}\sum_{x_0\in n^{-1}\Z^d,i}(H_{x_0,i})^2\Big)^{1/2}
\le c_3(\eta_1+\eps),\nonumber
\end{eqnarray}
if $\delta $ is taken small enough and $n$ is large. 
We thus complete the proof of (\ref{t7C}).\qed
  
\medskip

When $d=1$, Lemma \ref{cltlem} can be proved under much milder conditions.

\medskip
\noindent (A6) There exists $R>0$ and a Borel measurable $a:\R^d\to\cal M$ such that 
for each $r>0$
\begin{equation}\label{intcov}
\lim_{n\to\infty}\int_{|x|\le r}|a^n(x)-a(x)|dx=0.
\end{equation} 

\begin{cor}\label{d1cor}
Let $d=1$ and suppose (A1)-(A3) and (A6) hold. Then the conclusions of
Theorem \ref{clt} hold.
\end{cor}
\noindent {\sc Proof:} The proof is  similar to the proof of Theorem \ref{clt}. Let us point out 
the places where we need modifications. 
First, we can prove that there exist $c_1,c_2>0$ such that 
$c_1\le a^n(x)\le c_2$ for all $x\in \R^d$ and $n\in \N$. Indeed, by (A2) the lower bound
is guaranteed and the upper bound can be proved similarly to (\ref{aolacop}).
So, we know $\ce_{a^n}(f,f)$ is 
bounded whenever $f\in L^2$. For the proof that the right hand side of
(\ref{clt4}) goes to $0$ as $n\to\infty$, we use (\ref{intcov}).
(To be more precise, the  convergence of $a^n$ to $a$ locally in $L^2$ is used there, 
which is guaranteed by (\ref{intcov}) and the fact that the $a^n$ 
are uniformly bounded.) 
Noting these facts, the proofs of Theorem \ref{clt} and 
Proposition \ref{clttight} go the same way as above. For the proof of Lemma \ref{cltlem},
in Step 1, we do not need to control the oscillation of $a$ on each $\sS_m$.
Step 5 is not needed. We have 
that the expression (\ref{t72}) is equal to $\ol a^n_{ij_0}(x_0)$, and this is equal
to the expression in (\ref{t73q}). (This is a key point; because of this
we do not have to worry about the oscillation of $a$ and $a^n$.)  
Finally, in the computation of (\ref{compcdis}), the difference on the set
$B$ is $0$ due to the fact just mentioned, 
and we can prove that (\ref{compcdis}) is small directly. \qed

\medskip
We now give an extension of the result in \cite{SZ} to the case of
unbounded range. Assume

\medskip
\noindent (A7) \quad There exists $R>0$ such that for each $r>1$
\begin{equation}\label{szstrassmp}
\lim_{n\to\infty}\sum_{k\in \Z^d}\sup_{|y|\le nr}
\sup_{|x-y|\le nR}
\Big|C^{n,R}_{x,x+k}-C^{n,R}_{y,y+k}\Big|=0.
\end{equation}
\medskip
Let the $(i,j)$-th element of $b^n$  be given by
\begin{equation}\label{defsza}
\big(b^n(x)\big)_{ij}
=\sum_{k\in n^{-1}\Z^d} C_{nx, n(x+k)}^{n,R}n^2 k_ik_{j},\qquad x\in n^{-1}\Z^d.
\end{equation}
For general $x=(x_i)_{i=1}^d\in\R^d$, define $b^n(x):=b^n([x]_n)$.
Assume the $b^n$ version of (A6);

\medskip
\noindent (A8) There exists $R>0$ and a Borel measurable $a:\R^d\to\cal M$ 
such that for each $r>0$
\begin{equation}\label{intcovw2}
\lim_{n\to\infty}\int_{|x|\le r}|b^n(x)-a(x)|dx=0.
\end{equation} 
\medskip

We can recover and generalize the convergence theorem given in \cite{SZ} as follows. 
\begin{cor}\label{d1corsz}
Suppose that (A1)-(A3), (A7), and (A8)  hold. 
Then the conclusions of Theorem \ref{clt} hold.
\end{cor}
\noindent {\sc Proof:} 
For each $\eps>0$, let $R'=R'(\eps)>0$ be an integer that satisfies 
$\sum_{s\ge R'}\vp (s)s^2<\eps$.  Note that $C^{n,R}_{x,y}=
C^{n,R'/n}_{x,y}+1_{\{|x-y|>R'\}}C^{n,R}_{x,y}$. 
Then, 
for any $r\ge 1$,  
any $x\in n^{-1}\Z^d$ such that $|x|\le r$, and any $n\ge R'/R$,  
we have 
\begin{eqnarray*}
\Big|\big(a^n(x)\big)_{ij}-\big(b^n(x)\big)_{ij}\Big|
&\le &\sum_{k'\in \Z^d} 
\Big|\sum_{y:(y,k')\in L_x^{i,*}}C_{ny, ny+k'}^{n,R}\sgn k'_i-C_{nx, nx+k'}^{n,R}k'_i\Big|k'_j\\
&\le &R'^2\Big(\sum_{k'\in \Z^d}\sup_{|y'|\le nr}\sup_{|x'-y'|\le R'} 
\Big|C^{n,R'/n}_{x',x'+k'}-C^{n,R'/n}_{y',y'+k'}\Big|\Big)
+2\sum_{s\ge R'}\vp (s)s^2\\
&\le &R'^2\Big(\sum_{k'\in \Z^d}\sup_{|y'|\le nr}\sup_{|x'-y'|\le nR} 
\Big|C^{n,R}_{x',x'+k'}-C^{n,R}_{y',y'+k'}\Big|\Big)
+2\eps,
\end{eqnarray*}
where 
$L^{i,*}_z=\{(y,k')\in (n^{-1}\Z^d)\times \Z^d: 
y+n^{-1}\sP(k') \mbox{ contains the line segment from $z$ to } z+n^{-1}e^i\}$.
In the second inequality, we used the fact that if $(y,k')\in L^{i,*}_x$ and $x'=nx, y'=ny$, 
then $|x'-y'|=n|x-y|\le n|k'/n|=k'\le 
n\cdot R'/n=R'$.
Using (\ref{szstrassmp}) in (A7), 
the right hand side converges to $0$ as $n\to\infty$. 
In other words, 
\begin{equation}\label{aasznear}
|(a^n(x))_{ij}-b^n(x))_{ij}|\to 0~~\mbox{uniformly on compacts as } n\to\infty.
\end{equation}
Similarly, for any $r\ge 1$, we can prove 
\begin{equation}\label{aasz0}
|(b^n(x))_{ij}-(b^n(y))_{ij}|\to 0~~\mbox{ as }~~n\to\infty,~~
|x-y|\le n^{-1}R,~ |x|\le r.
\end{equation}
Now the proof of this corollary goes similarly to the proofs above. As before we point out places where we need modifications. 
First, as in Corollary \ref{d1cor}, we can prove that there exist $c_1,c_2>0$ such that 
$c_1I\le b^n(x)\le c_2I$ for all $x\in \R^d$ and $n\in \N$. So we know $\ce_{b^n}(f,f)$ is 
bounded whenever $f\in L^2$. As in Corollary \ref{d1cor}, we use 
(\ref{intcovw2})
to show that the right hand side of (\ref{clt4}) goes to $0$ as $n\to\infty$. 
Noting these facts, the proofs of Theorem \ref{clt} and 
Proposition \ref{clttight} go in the same way as before. For the proof of Lemma \ref{cltlem},
in Step 1, we do not need to control the oscillation of $a$ on each $\sS_m$.
Step 4 with respect to  $b^n$ works due to (\ref{aasznear}). Step 5 is not needed. 
Thanks to (\ref{aasznear}) and (\ref{aasz0}), the difference between 
the expression in (\ref{t72}) (with $a$ replaced by $b$) and the expression in (\ref{t73q}) is small. (This is again the key point; because of 
this we do not have to worry about the oscillation of $a$ and $b^n$.)  
Finally, in the computation of (\ref{compcdis}), the difference on the set
$B$ is small due to the fact just mentioned.\qed

\begin{rem}\label{remark2} {\rm If (A7) does not hold,  $b^n$ need not 
be the right approximation in general. 
Indeed, here is an example where $a^n$ converges to $a$, but $b^n$ does not as $n\to\infty$. 
Suppose $d=1$ and let $C^n_{k, k+i}$ equal $r_i$ if $k$ is odd, $s_i$ if 
$k$ is even, $i=1,2$. Then, we have 
\begin{eqnarray*}
b^n(k/n)
&=&\left\{\begin {array}{ll} r_1+s_1+8r_2,\qquad \mbox{if $k$ is odd},\\ 
r_1+s_1+8s_2,\qquad \mbox{if $k$ is even}.\end{array}\right.\\
a^n(k/n)
&=&\left\{\begin {array}{ll} 2r_1+4(r_2+s_2),\qquad \mbox{if $k$ is odd},\\ 
2s_1+4(r_2+s_2),\qquad \mbox{if $k$ is even}.\end{array}\right.\end{eqnarray*}
Suppose $r_1=s_1$  and $r_2\ne s_2$.
Then, the value of $b^n(k/n)$ depends on whether $k$ is odd or even,
so $b^n$ does not converge locally in $L^2$ as $n\to\infty$, whereas  
$a^n(k/n)=2r_1+4(r_2+s_2)$ is constant. 
In this case, the assumption of Theorem \ref{clt} (and Corollary \ref{d1cor})
holds and $a(x)=2r_1+4(r_2+s_2)$.
}
\end{rem}

\medskip

Theorem \ref{clt} gives a central limit theorem for the processes $Y^{(n)}$. 
Note that the base measure for $Y^{(n)}$ is the uniform measure, which 
converges with respect to
Lebesgue measure on $\R^d$. We finally discuss the convergence of 
the discrete time Markov chains $X^{(n)}$.
Let $Y^\nu_t$ be the continuous time $\nu$-symmetric Markov chain on $\Z^d$
which corresponds to $(\ce,\cf)$. It is a time change of  $Y_t$ and
it can be defined from $X_n$ as follows.
Let $\{U_i: i\in \N, x\in \Z^d\}$ be an independent collection of exponential random 
variables with parameter 1 that are independent of $X_n$.
Define $T_0=0, T_n=\sum_{k=1}^n U_k$. 
Set $\widetilde Y^\nu=X_n$ if $T_n\le t<T_{n+1}$; then 
the laws of $\widetilde Y^\nu$ and $Y^\nu$ are the same.
Let $\nu^D$ be a measure on $\cs$ defined by 
$\nu^D(A)=D^{-d}\nu (DA)$ for $A\subset \cs$. Since $\cs\subset \R^d$, 
we will regard $\nu^D$ as a measure on $\R^d$ from time to time.
By (A1), we see that 
$c_1\mu^D(A)\le \nu^D(A)\le c_2\mu^D(A)$ for all $A\subset \cs$ and all
$d$. So $\{\nu^D\}_D$ is tight and there is a convergent subsequence. We assume 
the following.

\medskip
\noindent (A9) There exists a Borel measure ${\bar \nu}$ on $\R^d$ such that
$\nu^D$ converges weakly to ${\bar \nu}$ as $D\to\infty$.

\medskip
Let $Z^{\bar \nu}_t$ be the diffusion process corresponding to the Dirichlet form
$\ce_a$ considered on $L^2(\R^d,{\bar \nu})$. It is a time changed process of
$Z_t$ in Theorem \ref{clt}. Note that by (A1), ${\bar \nu}$ is mutually absolutely
continuous with respect to Lebesgue measure on $\R^d$ so it charges no set of zero capacity.
Further, the heat kernel for $Z^{\bar \nu}_t$ still enjoys the estimates (\ref{aroest}).
 
Now we have a corresponding theorem for the discrete time Markov chains $X^{(n)}$.
Define 
\[ W^{(n)}_t=X^{(n)}_{[nt]}/\sqrt n. \]

\begin{cor}\label{cltcor}
Suppose (A1)-(A3),  (A5), and (A9) hold. 

(a) Then for each $x$ and each $t_0$ the  
$\P^{[x]_n}$-law of $\{W^{(n)}_t; 0\leq t\leq t_0\}$ converges weakly with respect to the topology of
the space $D([0,t_0], \R^d)$. The limit probability gives full
measure to  $C([0,t_0], \R^d)$.

(b) If $Z^{\bar \nu}_t$  is the canonical process on $C([0,\infty), \R^d)$ and $\P^x$ is the
weak limit of the $\P^{[x]_n}$-laws of $W^{(n)}$, then
the process $\{Z^{\bar \nu}_t, \P^x\}$ has continuous paths and is the symmetric process
corresponding  to the Dirichlet form $\ce_a$ considered on $L^2(\R^d,{\bar \nu})$. 
\end{cor}

\noindent {\sc Proof:} 
Let $Y^{(n),\nu}_t$ be the continuous time Markov chains on $\Z^d$ corresponding to
$\ce_n$ considered on $L^2(\Z^d,\nu)$, and set 
$Z^{(n),\nu}_t=Y^{(n),\nu}_{nt}/\sqrt n$. Then, by changing the measure $\mu^D$ to 
$\nu^D$ in the proof, we have the results corresponding to Theorem \ref{clt}
for $Z^{(n),\nu}_t$ and $Z^{\bar \nu}_t$. So it suffices to show
that there is a metric for $D([0,t_0], \R^d)$ with respect to which
the distance between 
$W^{(n)}$ and $Z^{(n),\nu}$ 
goes to 0 in probability, where in the definition of $Z^{(n),\nu}$ we use the realization
of $Y^{(n),\nu}$ given in terms of the $X^{(n)}$ by means of independent exponential random variables
of parameter 1. 

We use the $J_1$ topology of Skorokhod; see \cite{Bi}. The
paths of $Y^{(n),\nu}$ agree with those of $X^{(n)}$ except that the times of the jumps 
do not agree. 
Note that $X^{(n)}$ jumps at times $k/n$, while $Y^{(n),\nu}$ jumps at times
$T_k/n$.
So it suffices to show that if $T_k$ is the sum of i.i.d.\  exponentials with
parameter 1, then for each $\eta>0$ and each $t_0$
\[ \P(\sup_{k\leq [nt_0]} |T_k-k|\geq n\eta)\to 0 \]
as $n\to \infty$.
But by Doob's inequality, the above probability is bounded by
\[ \frac{4\Var T_{[nt_0]}}{n^2\eta^2}=\frac{4[nt_0]}{n^2\eta^2}\to 0 \]
as desired.
\qed

\begin{rem}\label{remark3}
{\rm We remark that the definition of $a^n$, and hence the statement of (A5), depends
on the definition of $\sP(k)$ and of the extension operator $E_n$. It would be nice
to have a central limit theorem with a more robust statement.
}
\end{rem}
\medskip

\begin{rem}\label{remark4}
{\rm 
We make a few comments comparing the central limit theorem in our
paper and the convergence theorem in \cite{SZ} in the case of 
bounded range. The result in \cite{SZ}
requires a smoothness condition on the conductances $C^n_{xy}$, while
we require smoothness instead on the $a^n$. Thus our theorem has
weaker hypotheses, and as Remark \ref{remark2} shows, there are 
examples where one set of hypotheses holds and the other set does not.
On the other hand, if 
(A1)-(A3) hold, then the $\{b^n\}$
will automatically
be symmetric, equi-bounded and equi-uniformly elliptic;
if in addition $b^n\to a$, then $a$ will be bounded
and uniformly elliptic  and this does not need to be assumed.
}
\end{rem}

\begin{flushright}
$\begin{array}{l}
\mbox{Richard F. Bass}\\
\mbox{Department of Mathematics}\\
\mbox{University of Connecticut}\\
\mbox{Storrs, CT 06269, U.S.A.}\\
\mbox{E-mail: {\tt bass@math.uconn.edu}}\\
\mbox{\ }\\
\mbox{Takashi Kumagai}\\
\mbox{Research Institute for Mathematical Sciences}\\
\mbox{Kyoto University}\\
\mbox{Kyoto 606-8502, Japan}\\
\mbox{E-mail: {\tt kumagai@kurims.kyoto-u.ac.jp}}\\
\end{array}$
\end{flushright}

\end{document}